\documentclass[11pt,a4paper]{article}
\usepackage[utf8]{inputenc}
\usepackage[english]{babel}
\usepackage{xcolor}
\usepackage{amsmath,amssymb,amsthm}
\usepackage[normalem]{ulem}

\theoremstyle{plain}
\newtheorem{theorem}{Theorem}[section]
\newtheorem{lemma}[theorem]{Lemma}
\newtheorem{corollary}[theorem]{Corollary}
\newtheorem{prop}[theorem]{Proposition}

\newtheorem{proposition}[theorem]{Proposition}
\newtheorem{example}[theorem]{Example}
\newtheorem*{theorema}{Theorem A}
\newtheorem*{theoremb}{Theorem B}
\newtheorem*{theoremc}{Theorem C}

\theoremstyle{definition}
\newtheorem{definition}[theorem]{Definition}
\newtheorem{remark}[theorem]{Remark}

\newcommand{\Cf}{\ensuremath \mbox{\large${\mathbf{1}}$}}

\newcommand{\G}{\ensuremath \mbox{\tiny $G$}}	
\newcommand{\dmg}{\ensuremath  \, \mathrm{dm}_{\G}}
\newcommand{\Mg}[1]{M_{\mbox{\tiny $\widehat{#1}$\,}}}

\newcommand\<[1]{\left\langle\, #1\,\right\rangle}
\newcommand{\N}{\mathbb{N}} 
\newcommand{\Z}{\mathbb{Z}} 
\newcommand{\Q}{\mathbb{Q}} 
\newcommand{\R}{\mathbb{R}} 
\newcommand{\C}{\mathbb{C}} 
\newcommand{\T}{\mathbb{T}}

\newcommand\norm[1]{\lVert #1 \rVert}

\newcommand{\mA}{\ensuremath  \mathfrak{A}}
\newcommand\supp[1]{\mathrm{supp}\left(#1 \right)}

\renewcommand{\leq}{\leqslant}
\renewcommand{\geq}{\geqslant}

\newcounter{equi1}

\newcommand\restr[2]{{
		\left.\kern-\nulldelimiterspace 
		#1 
		\vphantom{|} 
		\right|_{_{#2}} 
}}

\title{Recurrent and rigid multipliers on commutative Banach algebras.  The Fourier algebra}

\author{María J. Beltrán-Meneu\footnote{Instituto Universitario de Matem\'aticas y Aplicaciones (IMAC), Universitat Jaume I,
Av. Vicent Sos Baynat, s/n, Castellón de la Plana, E-12071 (Spain), mmeneu@uji.es,   https://orcid.org/0000-0001-9645-096X.},
	Jorge Galindo\footnote{Instituto Universitario de Matem\'aticas y Aplicaciones (IMAC), Universitat Jaume I,
		Av. Vicent Sos Baynat, s/n, Castellón de la Plana, E-12071 (Spain), jgalindo@uji.es,   https://orcid.org/0000-0002-5680-1742.}
	and Enrique Jordá\footnote{Instituto Universitario de Matem\'atica Pura y Aplicada IUMPA, Universitat Politècnica de València, Camino de Vera, s/n, Valencia,  E-46022 (Spain), ejorda@mat.upv.es, https://orcid.org/0000-0003-2980-1699.}}

\date{}

\begin{document}
	
	\maketitle
	
	\begin{abstract}
%

    We study recurrence, rigidity and uniform rigidity of multipliers on semisimple commutative Banach algebras. For a power-bounded multiplier $T$ on a Tauberian, self-adjoint algebra $\mA$, we show that  $T$ satisfies each of these properties if and only if  the multiplication operator on $C_0(\Delta(\mA))$ induced by the Gelfand transform of $T$  does.  One consequence of our analysis is that, if $\Delta(\mA)$ is $\sigma$-compact, then recurrence implies rigidity for such multipliers.  We examine in detail the multipliers of the Fourier algebra $A(G)$. In particular, we show that recurrence of a multiplier of $A(G)$ is equivalent to recurrence of the corresponding multiplier on $C_0(G)$ when the connected components of $G$ are open. For multipliers that are not power bounded,  a flexible construction based on chains of compact subgroups of a locally compact abelian group $G$ produces  multipliers of $A(G)$ that are rigid but not uniformly rigid, and multipliers that are recurrent but not rigid, such that in each case the corresponding multiplier on $C_0(G)$ is uniformly rigid. Via the Fourier-Stieltjes transform, these results describe the recurrence and rigidity properties of convolution operators on $L^1(G)$.

	\end{abstract}
	
	

\section{Introduction}

The study of recurrence in dynamical systems dates back to the seminal work of Poincaré and Birkhoff at the end of the 19th century and early 20th century. In its classical formulation, a transformation is said to be recurrent if almost all points in the space eventually return to positions arbitrarily close to their initial one. Its systematic exploration within the framework of linear dynamics is a relatively recent development, initiated by Costakis, Manoussos, and Parissis in their  papers
\cite{costamanopari12,costamanopari14}, with roots in the papers of Eisner and Grivaux \cite{eisn09,eisngriv11}.

A bounded linear operator $T$ acting on a Banach space $X$ is said to be \textit{recurrent} if, for every non-empty open set $U \subseteq X$, there exists an integer $k \in \mathbb{N}$ such that $U \cap T^{-k}(U) \neq \emptyset$. Equivalently, $T$ is recurrent if the set of its recurrent vectors, that is, those $x \in X$ for which there exists a sequence $(k_n)_n$ of positive integers such that $T^{k_n}x \to x$ as $n \to \infty$, is dense in $X$. If $T$ is power-bounded, the set of recurrent vectors  is closed \cite[Lemma  3.1]{costamanopari14}, hence $T$  is recurrent if and only if every vector in $X$ is recurrent.

Recurrence provides a more flexible framework for  the dynamics of linear operators on infinite-dimensional spaces than the classical theory of hypercyclicity. This contrast is clearly illustrated by unitary operators on Hilbert spaces: although they are never hypercyclic, they  display a rich  recurrent behaviour.

Alongside recurrence, the concepts of rigidity and uniform rigidity have also received significant attention. A bounded linear operator $T$ is said to be \emph{rigid} if there exists an increasing sequence $(k_n)_n$ of positive integers such that $T^{k_n} \to I$ in the strong operator topology, and \emph{uniformly rigid} if this convergence holds in the operator norm. In either case, we say that $T$ is rigid (respectively, uniformly rigid) \emph{along the sequence} $(k_n)_n$. As shown in \cite{costamanopari14}, these properties  capture subtle differences in the asymptotic behavior of linear orbits.

The present paper addresses the problem of recurrence and rigidity for multipliers of semisimple commutative Banach algebras (see  \cite{bulancea25} for initial explorations of  this theme in the unital setting). Our approach has two main objectives.  First, we relate the recurrence and rigidity properties of a multiplier $T\colon \mA\to \mA$ on a commutative Banach algebra $\mA$  to those of  the multiplication operator $\Mg{T}\colon C_0(\Delta(\mA))\to C_0(\Delta(\mA))$ that has as symbol  its Gelfand transform $\widehat{T}$. Together with our analysis of the particular case  of the algebras $C_0(\Omega)$, where $\Omega$ is a locally compact Hausdorff space,  this yields clean and effective characterizations of recurrence and rigidity of multipliers, at least when $T$ is power-bounded and $\mA$ is self-adjoint and Tauberian. Second, we use these characterizations to classify the recurrence and rigidity of the  multipliers of the convolution Banach algebras $L^1(G)$, where $G$ is a locally compact Abelian group, and their analogue for noncommutative groups, the Fourier algebra $A(G)$, where $G$ is a locally compact group.

We now state our main results; see Section \ref{term} for the mostly standard terminology used here.

Theorem A will be proved in Lemma \ref{recvec} and Theorems  \ref{recA}, \ref{rigidA}, \ref{sigmacompact},  and \ref{URA}.

\begin{theorema}\label{tha}
Let $\mA$ be a semisimple and self-adjoint commutative  Banach algebra with spectrum $\Delta(\mA)$ and let $ T\colon \mA\to \mA$ be a power bounded multiplier with Gelfand transform $\widehat{T}$.
\begin{enumerate}
\item Assume that  $\widehat{T}(\Delta(\mA))\subseteq\T$. For $a\in\mA$, set $C=\supp{\widehat{a}}\subseteq\Delta(\mA)$. If $\widehat{T}(C)$ is a Dirichlet subset of $\T$, then $a$ is a recurrent vector for $T$.
\item If $\mA$ is Tauberian, then $T$ is recurrent (respectively, rigid or uniformly rigid) if and only if  the operator $\Mg{T}\colon C_0(\Delta(\mA))\to C_0(\Delta(\mA))$ given by $\Mg{T}f=\widehat{T}\cdot f$, is recurrent (respectively, rigid or uniformly rigid).
\item If $\mA$ is Tauberian  and $\Delta(\mA)$ is $\sigma$-compact, then $T$ is recurrent if and only if it is rigid.
\end{enumerate}
\end{theorema}

The relevance of this theorem lies in the fact that it draws consequences in the $\mA$-norm from data in the spectrum. When power boundedness, the Tauberian property, or $\sigma$-compactness are not assumed in the respective statements, their conclusions need not  hold. We will provide  examples to this effect.

Assertion~3 gives a partial positive answer to \cite[Question~9.5]{costamanopari14}, which, in view of \cite[Proposition~3.2(ii)]{costamanopari14}, can be stated as follows: Is every power-bounded recurrent operator on a separable Banach space rigid? Our result settles the question affirmatively for multipliers of Tauberian, self-adjoint Banach algebras.

A multiplier $T$ of $A(G)$, which we denote by $T\in M(A(G))$, is always a multiplication operator $M_\psi$, $u\mapsto \psi\cdot u$, induced by a
continuous function $\psi\colon G\to \C$. Of special interest are the multipliers $M_\phi$ with $\phi\in B(G)$, the Fourier-Stieltjes algebra of $G$; if $G$ is amenable there are actually no other multipliers.  As $\Delta(A(G)) = G$,  both $T=M_\psi$ and $M_{\widehat{T}}$ are multiplication operators with the same symbol $\psi$. We prefer to denote  $\Mg{T}$ as $M_{\psi}^0 \colon C_0(G)\to C_0( G)$.

Theorem B summarizes our results for the Fourier algebra $A(G)$ of a locally compact group $G$; its proof is spread over the results of Section 5.
\begin{theoremb}\label{thb}
Let  $G$ be a locally compact group and let $T=M_\psi$ be a multiplier of $A(G)$.
\begin{enumerate}
\item If the connected components of $G$ are open, then $M_\psi$ is recurrent if and only if $\psi(G)\subseteq \T$ and $\psi$ is constant on each connected component.
\item If $M_\psi$ is power bounded, then $M_\psi$ is recurrent (respectively rigid, uniformly rigid) if and only if $M_\psi^0$ is so.
\item If $G$ is $\sigma$-compact and $M_\psi$ is power bounded, then $M_\psi$ is recurrent if and only if it is rigid.
\item  If $\phi\in B(G) $ is power bounded and $M_\phi$ is recurrent, then $\phi$ is a piecewise affine function $\phi \colon G\to  \T$;
     if $\|M_{\phi}\|_{M(A(G))}\leq 1,$  then $\phi$ must even  be affine.
\item If $\phi\in B(G) $ is power bounded, $M_\phi$ is uniformly rigid if and only if $\phi(G)$ is a finite subset of $\T$. This need  not hold  when $\phi$ is only power-bounded in the multiplier norm.
\item If $\phi\in B(G)$ and $M_\phi$ is not power bounded, the recurrence, rigidity, and uniform rigidity properties of $M_\phi$ and $M_\phi^0$ can be different. For the rigidity and uniform rigidity cases, we obtain examples on countable  groups.    Also, $M_\phi$ can be recurrent but fail to be rigid, in contrast with 3.
\end{enumerate}
\end{theoremb}

We note that the examples in assertion 6 of Theorem B all follow from a single construction that uses  sums of alternating blocks of Haar measures of a chain of compact subgroups. Depending on the length of these sums,  the Fourier-Stieltjes transform produces different symbols with the same infinite Dirichlet range but different dynamical properties: they can  yield uniformly rigid multipliers,  rigid but not uniformly rigid multipliers and recurrent but  non-rigid multipliers  (see Proposition~\ref{abelian-block-construction}).

When $G$ is a locally compact Abelian group, the Fourier–Stieltjes transform identifies $L^1(G)$ with $A(\widehat{G})$ and the multipliers of $L^1(G)$   (i.e., convolution operators by bounded measures in $M(G)$) with $B(\widehat{G})$. The conclusions of Theorem B can be stated in terms  of convolution operators on the group algebra $L^1(G)$. We state them here in the particular case of $L^1(\T)$.

\begin{theoremc}
	Let $\mu\in M(\mathbb{T})$ and consider the convolution operator $C_\mu\colon L^1(\mathbb{T})\to L^1(\mathbb{T}),$  $C_\mu f=\mu*f$.
	\begin{enumerate}
		\item[(i)] $C_\mu$ is recurrent if and only if $|\widehat\mu(k)|=1$ for every $k\in\Z$.
		\item[(ii)] If $C_\mu$ is power bounded, it is recurrent if and only if it is rigid, and it is uniformly rigid if and only if $\widehat\mu(\Z)$ is a finite subset of $\T$.
		\item[(iii)] There exist non-power-bounded convolution operators on $L^1(\T)$ which are uniformly rigid, rigid but not uniformly rigid, and recurrent but not rigid. All three can be chosen with the same Fourier--Stieltjes range $\widehat\mu(\Z)=\{e^{2\pi i/n!}:n\geq1\}$; their associated operators on $c_0(\Z)$ are uniformly rigid.
	\end{enumerate}
\end{theoremc}

\section{Preliminaries}\label{term}
	\subsection{Sparse subsets of $\T$}
		
	As established in the introduction, recurrence and (uniform) rigidity of multipliers on commutative Banach algebras are closely connected to the geometry of subsets of $\mathbb{T}$, particularly the behaviour of their powers. Throughout the paper, we will make repeated use of Dirichlet sets, already introduced in the introduction, and will also find it useful to consider Arbault sets. We now formally define these notions.
	
	\begin{definition}
		Let $E$ be a subset of the unit circle $\T$.
		\begin{itemize}
			\item $E$ is a \emph{Dirichlet} set   if there is an increasing sequence $(n_k)_k\subseteq \N$ such that 	\\		$\lim_k \sup_{x\in E} \left| x^{n_k}-1\right|=0$.
			\item $E$ is an \emph{Arbault} set if there is an increasing sequence $(n_k)_k\subseteq \N$ such that, for any $x\in E$, 	$ \lim_k \left| x^{n_k}-1\right|=0.$
		\end{itemize}
	\end{definition}
	Finite sets are Dirichlet; this is actually \emph{Dirichlet's theorem} \cite[Chapter I]{lind69},  whereas every countable subset of $\T$ is known to be an Arbault set \cite{dikrkunen07}.  The term Arbault set stems from the paper \cite{arba52} where they were first studied;     see also \cite{barbetal17} and the references therein. There are Arbault sets that are not Dirichlet, any countable dense subset of $\T$ is such an example. A classical example of an infinite Dirichlet set is $\{e^{\frac{2\pi i}{n!}}\colon n\in \N\}$, as witnessed by the sequence $n_k=k!$.

	\subsection{Multipliers of  commutative Banach algebras}

We will be working with a commutative Banach algebra  $\mA$.   The set of all   non-zero algebra homomorphisms from $\mA$ onto $\C,$ which is a locally compact Hausdorff space (compact, if $\mA$ is unital) will be denoted by $\Delta(\mA)$.
 The space $\Delta(\mA)$ is equipped with the $w^*$-topology of $\mA^*$ and is called the \emph{spectrum} or \emph{maximal ideal space} of $\mA.$

For $a \in \mA$, the continuous function $\hat{a} : \Delta(\mA) \to \mathbb{C}$, given by  $\hat{a}(\varphi) = \varphi(a)$ is called the \emph{Gelfand transform} of $a$. If $C_0(\Delta(\mA))$ is considered with its usual (pointwise multiplication, sup-norm) Banach algebra structure, the correspondence $a\mapsto \Gamma(a)=\widehat{a}$  defines  an algebra homomorphism  $\Gamma : \mA \to C_0(\Delta(\mA))$, the \emph{Gelfand representation} of $\mA$. Since   $\|\hat{a}\|_\infty = r(a)\leqslant \|a\|$ for every $a \in \mA$, $\Gamma$ is a  continuous linear map. It is isometric if and only if $\|a^2\| = \|a\|^2$ for all $a \in \mA$. In this case, we say the algebra $\mA$ is \emph{uniform.}

We say that $\mA$ is \emph{semisimple} if, for every nonzero $a\in\mA$, there exists $\chi\in\Delta(\mA)$ such that $\widehat{a}(\chi)=\chi(a)\neq 0$. Equivalently, the Gelfand representation $a\mapsto\widehat{a}$ is injective.

The algebra $\mA$ is said to be \emph{regular} if given a closed subset $E$ of $\Delta(\mA)$ and some $\chi\in \Delta(\mA)\setminus E,$ there exists $a\in \mA$ such that $\widehat{a}=0$ on $E$ and $\widehat{a}(\chi)\neq 0.$ Every regular commutative Banach algebra is \emph{normal} in the sense that whenever $E\subseteq \Delta(\mA)$ is closed, $K\subseteq\Delta(\mA)$ is compact and $E\cap K=\emptyset,$ then there exists $a\in \mA$ such that $\supp{\widehat{a}}\subseteq \Delta(\mA)\setminus E$ and $\widehat{a}|_K=1.$

When  the set of all $a\in \mA$ such that $\widehat{a}$ has compact support is dense in $\mA$ we say that $\mA$ is  \emph{Tauberian}.

   $\mA$  is said to be  \emph{self-adjoint} if for every $a\in \mA$ there exists $b\in \mA$ such that $\hat{b}=\overline{\widehat{a}},$ where $\overline{\widehat{a}}$ denotes the complex conjugate of $\widehat{a}\in C_0(\Delta(\mA)),$ that is,   $\overline{\widehat{a}}(\chi) = \overline{\widehat{a}(\chi)}   \text{ for all } \chi \in \Delta(\mA).$
 If $\mA$ is semisimple, this element is unique. In this case we denote it by $\overline{a}$ and we have $\widehat{\overline{a}}=\overline{\widehat{a}}.$ A semisimple and  self-adjoint Banach algebra $\mA$  is a Banach *-algebra with   conjugation as involution, that is, $a^*:=\overline{a},$ $a\in \mA.$ Indeed, there exists $C>0$ such that
 \begin{equation} \label{inv}
 	\frac{1}{C}\|a\|\leq  \|\overline{a}\|\leq C \|a\|
 \end{equation}

\noindent  for every $a\in \mA$ \cite[Corollary 2.1.12]{kaniuth_CommBA}.


A linear operator $T:\mA\rightarrow \mA$  is a \emph{multiplier} of $\mA$ when  $T(ab)=aT(b)$ for all $a,b\in \mA.$

In what follows, we assume that $\mA$ is a semisimple commutative Banach algebra. In this case, every multiplier of $\mA$ is bounded, and the set $M(\mA)$ of all multipliers of $\mA$, called the \emph{multiplier algebra} of $\mA$, is a unital commutative Banach subalgebra of $\mathcal{B}(\mA)$, the Banach algebra of bounded linear operators on $\mA$.

For every $a\in \mA$, the operator $T_a:\mA\rightarrow \mA,$ $b\mapsto a\cdot b$ is a multiplier of $\mA$. Since the algebra is semisimple,  $\mA$ can be identified with a subalgebra of $M(\mA)$.  If $\mA$ is unital, we even get $M(\mA)=\mA$, as for every $T\in M(\mA),$  $T=T_a$, with $a=T(1).$

The Gelfand representation of $\mA$ can be extended to a continuous algebra homomorphism
\begin{align*}
  \Gamma_M& : M(\mA) \to C_b(\Delta(\mA)),
\end{align*} where, for a given $T\in M(\mA)$, $\Gamma_M(T)$ is usually denoted by $\widehat{T}$ and defined as the  uniquely determined  continuous function satisfying 	$$\widehat{T(a)}=\widehat{T}\cdot \widehat{a}$$
for all $a\in \mA$.   Observe that if $\mA$ is semisimple, then $\Gamma_M$ is injective, and evaluating $\Gamma_M$ at the translation operators $T_a$ shows that $\restr{\Gamma_M}{\mA} : \mA \to C_0(\Delta(\mA))$ is simply  the Gelfand representation  $\Gamma.$

With the aid of this representation, every multiplier $T$  can be associated to a multiplier of $C_0(\Delta(\mA)) $  \[\Mg{T}:C_0(\Delta(\mA))\rightarrow C_0(\Delta(\mA))\] given  by  $\Mg{T}(f)=\widehat{T}\cdot f$.



 \begin{prop}\label{Remark_inversa}
	Let  $T:\mA\rightarrow \mA$ be a multiplier such that  $\widehat{T}(\Delta(\mA))\subseteq \T$. Then  $T$ is injective and $\overline{T(\mA)}=\mA.$ If in addition $\mA$ is self-adjoint, then $T$ is invertible,  and $T^{-1}$ is power bounded if and only if $T$ is power bounded.
\end{prop}

\begin{proof}
	Let us see first that  $T$   is injective and $\overline{T(\mA)}=\mA.$ Indeed,  if $T(a)=0,$ then $\widehat{T}\cdot \widehat{a}=0,$ which, as $\mA$ is semisimple, yields $a=0.$ For the second statement, $T(\mA)$ is an ideal in $\mA$, if it is not dense, there exists $\chi \in \Delta(\mA)$ such that $\overline{T(\mA)} \subseteq \ker(\chi) \subseteq \mA$. Taking into account that $\widehat{T}(\chi) \widehat{u}(\chi) = \widehat{T(u)}(\chi) = \chi(T(u)) = 0$ for each $u \in \mA$, and that $\widehat{T}(\Delta(\mA)) \subseteq \mathbb{T}$, it is enough to choose $u \in \mA$ such that $\widehat{u}(\chi) \neq 0$ to reach a contradiction.
	
	Let us now check that $T$ is invertible by finding its inverse.
	Pick $a\in \mA$, then, for each $\chi\in \Delta(\mA)$,
	\begin{align*}
		\widehat{T\left(\overline{T(\overline{a})}\right)}(\chi)&=
		\widehat{T}(\chi)\widehat{\overline{T(\overline{a})}}(\chi)
		\\&=\widehat{T}(\chi)
		\overline{\widehat{T(\overline{a})}(\chi)}\\
		&=\widehat{T}(\chi)\overline{\widehat{T}(\chi) \widehat{\overline{a}}(\chi)}
		\\&=\widehat{a}(\chi).
	\end{align*}
	Since $\mA$ is semisimple and $T$ is  injective  we get $T^{-1}(a)=\overline{T(\overline{a})}$.
	%
	It defines a bounded operator on $\mA$ by (\ref{inv}). The power-boundedness assertion also follows from (\ref{inv}).
\end{proof}

%
%
%

For more details on commutative Banach algebras and multipliers on them, see e.g.  \cite{kaniuth_CommBA} and \cite{Larsen1971}.

We finish this section with   two preliminary results which connect the recurrence properties of a multiplier $T\in M(\mA)$ with the range of $\widehat{T}.$

\begin{proposition}
	\label{lem:T}
Let $\mA$ be a Banach algebra and let $T\in M(\mA)$. If $T$ is recurrent then $\widehat{T}(\Delta(\mA))\subseteq \T$.
\end{proposition}
\begin{proof}
	Let $\chi \in \Delta(\mA)$.  Since $\{a\in \mA:\ \chi(a)\neq 0\}$ is open and non-empty,   we can get  a recurrent vector $a\in \mA$ such that $\chi(a)\neq 0$. Now we have
	$$|(\widehat{T}(\chi)^n-1)\widehat{a}(\chi)|=|\widehat{(T^n(a)-a)}(\chi)|\leq \|\widehat{T^n(a)-a}\|_\infty\leq \|T^n(a)-a\|.$$

\noindent Recurrence of $a$ allows us to get an increasing sequence $(n_k)$ of natural numbers such that 			$|\widehat{T}(\chi)^{n_k}-1|$ converges to 0. This certainly implies $\widehat{T}(\chi)\in\T$.
\end{proof}

\begin{proposition}
	\label{rangofinito}
Let $\mA$ be a semisimple Banach algebra and let $T\in M(\mA)$ such that $\widehat{T}(\Delta(\mA))\subseteq \T$ is a finite set. Then $T$ is power bounded and uniformly rigid.
\end{proposition}
\begin{proof}
	
The orbit $\{\widehat{T}^n: n\in\N_0\}\subseteq \text{span}\{\Cf_{\widehat{T}^{-1}(\{z\})}: \ z\in \widehat{T} (\Delta(\mA)) \}$ spans a finite dimensional space in $C_b(\Delta(\mA))$, and, as $\Gamma_M$ is  injective,  also the span of $\{T^n: n\in\N_0\}$ is of  finite dimension in $M(\mA)$. Hence there exists $M>0$ such that, for each $n\in\N$ we have
	
	\begin{equation}\label{1}\|T^n\|\leq M \|\widehat{T}^n\|_\infty  \end{equation}
	 \noindent and also
	\begin{equation}\label{2}\|T^n-I\|\leq M \|\widehat{T}^n-1\|_\infty \end{equation}
		
	\noindent From \eqref{1} we get that $T$ is power bounded and from the hypothesis $\widehat{T}(\Delta(\mA))$ is finite in $\T$, Dirichlet theorem  and \eqref{2} we conclude that $T$ is uniformly rigid.
\end{proof}


\subsection{The group algebra and the Fourier algebra}

In our final section, we address the particular case of the convolution algebra $L^1(G)$ of a locally compact Abelian group $G$. As it is well known (and first proved by Wendel \cite{wend52}), the multipliers of this algebra correspond to convolution operators  $C_\mu\colon L^1(G)\to L^1(G)$  induced by bounded measures $\mu \in M(G)$ and defined by  $C_\mu(f)(t)= (\mu \ast f)(t)=\int_{G} f(x^{-1}t)\,\mathrm{d}\mu(x)$.

The usual description of  $\Delta(L^1(G))$ is based on the Fourier transform and requires  the duality of locally compact Abelian groups. We now briefly introduce both tools.

 If $G$ is a locally compact Abelian group,  its character group $\widehat{G}$ is the group of continuous homomorphisms of $G$ into the unit circle $\T$. Equipped with the topology of uniform convergence on compact subsets of $G$, $\widehat{G}$ is again  a locally compact group. Pontryagin's duality theorem establishes that $\widehat{\,\widehat{G}\;\,}$ is topologically isomorphic to $G$ in a canonical way.

If $f\in L^1(G)$, its  \emph{Fourier transform}  is a function $\widehat{f}\colon \widehat{G}\to \C$, defined as $\widehat{f}(\chi)=\int_{G} \chi(t^{-1})f(t) \, \dmg(t)$, where $\dmg$ is the Haar measure of $G$. This transform can  be used to obtain a realization of $\Delta(L^1(G))$: for each $\psi\in \Delta( L^1(G)), $ there is  $\chi\in \widehat{G}$ such that $\psi(f)=\widehat{f}(\chi)$, for every $f\in L^1(G)$. Identifying in this way $\Delta(L^1(G))$ with $\widehat{G}$, the Gelfand representation of $L^1(G)$ is just the Fourier transform, i.e., $\Gamma(f)=\widehat{f}$ for each $f\in L^1(G)$. This identification can be extended to the Gelfand representation $\Gamma_M$ of the  multiplier algebra $M(G)$ of $L^1(G)$. In this case,   $\Gamma_M(\mu)=\widehat{\mu}$, for every $\mu\in M(G)$, where $\widehat{\mu}$ is the \emph{Fourier-Stieltjes} transform, $\widehat{\mu}(\chi)=\int_{G} \chi(t^{-1}) \, \mathrm{d}\mu(t)$.

The range of the Fourier transform (i.e., $\Gamma(L^1(G))$)  is known as  the \emph{Fourier algebra} of $\widehat{G}$, denoted by $A(\widehat{G})$, while the range of the Fourier-Stieltjes transform (i.e., $\Gamma_M(M(G))$)  is known as  the \emph{Fourier-Stieltjes algebra} $B(\widehat{G})$. Via the Gelfand map $\Gamma_M$, we view $A(\widehat{G})$ and $B(\widehat{G})$ as Banach algebras isometric to $L^1(G)$ and $M(G)$, respectively. Moreover, since  the  character group of $\widehat{G}$
is topologically isomorphic to $G$, we can bypass duality and frame the multiplier theory of $L^1(G)$ (convolution operators by measures) as the multiplier theory of $A(G)$ (multiplication operators by functions in $B(G)$). We will adopt this perspective.

Working within $A(G)$, one need not restrict to commutative groups. Associated with any locally compact group, one can define two Banach function algebras on $G$ that coincide with $A(G)$ and $B(G)$ when the group $G$ is Abelian. These spaces are also termed the Fourier and Fourier-Stieltjes algebras of $G$; see \cite{eyma64} and Chapter 2 of \cite{kanilaubook} for their actual definitions and main properties. The spectrum of $A(G)$ coincides with point evaluations by members of $G$ and, since $A(G)$ is a closed ideal of $B(G)$, every element of $B(G)$ acts as a multiplier of $A(G)$. When $G$ is \emph{amenable}, it admits no other multipliers.   In nonamenable groups, however, $A(G)$ admits multipliers that do not originate from $B(G)$ (see \cite[Chapter 5]{kanilaubook}).
We will not dwell further here on the properties of these algebras---some other useful properties  are delayed to the introduction of Section~\ref{Section:A(G)}---but we can already indicate that our results on the dynamics of Fourier algebras  are valid in the context of locally compact noncommutative groups.

 \section{Multiplication operators on $C_0(\Omega)$}\label{Section:C_0}
One of our  aims is to compare the properties of a multiplier  $T$ as an operator on a Banach algebra $\mA$  with those of the multiplication operator $\Mg{T}$ acting on  $C_0(\Delta(\mA))$. For this to be meaningful, the rigidity  and recurrence properties of multiplication operators on $C_0(\Omega) $,  where $\Omega $ is  a locally compact Hausdorff space, should be understood. In this section we characterize these properties for the operator $M_{\phi}
\colon C_0(\Omega)\rightarrow C_0(\Omega),$ $f\mapsto \phi  \cdot f,$  where   $\phi\in C_b(\Omega)$ is a bounded continuous function on $\Omega$, that is, for all multipliers of $C_0(\Omega)$.  This characterization is linked to  the existence of Dirichlet sets within certain subsets of the range of $\phi$.

Since the evaluation functionals at points of $\Omega$  exhaust  $\Delta(C_0(\Omega))$,   we obtain  the
 following result as a corollary of Proposition \ref{lem:T}.
	
	\begin{proposition}\label{lem:T_C0}
Let $\Omega$ be a locally compact Hausdorff   space and  $\phi\in  C_b(\Omega)$.  If $M_{\phi}:C_0(\Omega)\rightarrow C_0(\Omega)$  is recurrent, then $\phi(\Omega)\subseteq \T$. Hence every recurrent multiplier on $C_0(\Omega)$ is an isometry.
	\end{proposition}

	\begin{prop}
		\label{nonunitalcarac}
Let $\Omega$ be a locally compact Hausdorff   space and  $\phi\in  C_b(\Omega)$ such that $\phi(\Omega)\subseteq \T$.
		\begin{itemize}
			\item[(i)]  $M_{\phi}:C_0(\Omega)\rightarrow C_0(\Omega)$  is recurrent if and only if $\phi(K)$ is a Dirichlet set for every compact set $K\subseteq \Omega$.
			\item[(ii)] $M_{\phi}:C_0(\Omega)\rightarrow C_0(\Omega)$  is  rigid if and only if there is an increasing sequence $(n_k)$ of natural numbers such that $\lim_k\sup_{x\in K} |\phi^{n_k}(x)-1|=0$  for every compact set $K\subseteq \Omega$.
			\item[(iii)] 	$M_{\phi}:C_0(\Omega)\rightarrow C_0(\Omega)$  is uniformly  rigid if and only if $\phi(\Omega)$ is Dirichlet.
		\end{itemize}
	\end{prop}
	\begin{proof}
(i) As $\phi(\Omega)\subseteq \T$, $M_{\phi}$ is power bounded.	As a consequence,  if  $M_{\phi}$ is recurrent, every $f\in C_0(\Omega)$ has to be a recurrent vector. Hence, for every compact $K\subseteq \Omega$, by applying the definition to any function $f$ such that $f|_K\equiv 1$, we get that $\phi(K)$ is Dirichlet.
		
Assume now that  $\phi(K)$ is Dirichlet for every $K\subseteq \Omega$ compact. Let $f\in C_0(\Omega)$ and $(K_n)_n$ an increasing sequence  of compact subsets of $\Omega$ satisfying $\sup_{x \in \Omega\setminus K_n}|f(x)|<\frac1n$, and $(m_n)_n$  an increasing sequence of natural numbers satisfying  $\sup_{x \in K_n}|{\phi}^{m_n}(x)-1|<\frac1n$. Then, as $\sup_{x\in \Omega\setminus K_n}|\phi^{m_n}(x)-1| |f(x)|\leq \frac{2}{n}$ and $\sup_{x\in  K_n}|\phi^{m_n}(x)-1| |f(x)|\leq \frac{\|f\|_{\infty}}{n}$, we get $\|\phi^{m_n}f-f\|_{\infty}\rightarrow 0$ as $n$ tends to infinity, as we wanted to see.
		
(ii) If $(m_n)_n$ is an increasing sequence of natural numbers such that $\lim_n M_{\phi}^{m_n}(f)=f$ for all $f\in C_0(\Omega)$, by applying this to functions identically one in  $K\subseteq \Omega$ compact, we get that $\phi(K)$ is Dirichlet along  $(m_n)_n$, and $(m_n)_n$ does not depend on $K$.  Conversely, if such a sequence   $(m_n)_n$ exists, for each $f\in C_0(\Omega)$ and $\varepsilon>0$ we choose $K_\varepsilon$ such that $\sup_{x\in \Omega\setminus K_\varepsilon}|f(x)|<\frac{\varepsilon}{2}$. Consider $n_0\in \N$ such that $\sup_{x \in K_\varepsilon}|\phi^{m_n}(x)-1|<\frac{\varepsilon}{ 2\|f\|_{\infty}}$ for  $n\geq n_0$. If $n\geq n_0$ then $\|M_{\phi}^{m_n}f-f\|<\varepsilon$. Therefore, $M_{\phi}$ is rigid.

		(iii)  As $\|M_{\phi}^n-I\|\leq \|\phi^n-1\|_{\infty},$ we get that $\phi(\Omega)$ Dirichlet  is a sufficient condition for the  uniform rigidity of $M_{\phi}$. To see that it is also necessary, consider  the sequence $(n_k)_k$   which gives  uniform rigidity. For a given $\epsilon >0$ there exists then $k_0\in \N$ such that $\|\phi^{n_k} f-f \|_{\infty}\leq  \epsilon$ for every $k\geq k_0$ and every $f\in C_0(\Omega),$ $\|f\|_{\infty}\leq 1$. For each  $x\in \Omega,$ take $f_x\in C_0(\Omega)$ such that $f_x(x)=1$ and $\|f_x\|_\infty \leq 1.$ Then,
\[|\phi^{n_k}(x)-1|=|\phi^{n_k}(x)f_x(x)-f_x(x)|\leq \|\phi^{n_k}f_x-f_x\|_\infty\leq \epsilon,\] which gives the conclusion.
	\end{proof}

When $\Omega$ is $\sigma$-compact, Proposition \ref{nonunitalcarac} yields  that some of  the dynamical properties we are studying are indistinguishable. Let us focus first on the compact case, where they all coincide. Here $C_0(\Omega)=C_b(\Omega)=C(\Omega).$

\begin{corollary}\label{equiv_compact}
Let $K$ be a   compact Hausdorff space and  $\phi\in  C(K)$.  Then the following statements are equivalent:
\begin{itemize}
	\item [(i)]	$M_{\phi}:C(K)\rightarrow C(K)$   is recurrent.
	\item [(ii)]	$M_{\phi}:C(K)\rightarrow C(K)$   is uniformly rigid.
		\item [(iii)] $\phi(K)$ is a Dirichlet subset of $\T.$
		\item[(iv)] There exists  $(n_k)_k$ such that $\lim_k\|\phi^{n_k}-1\|_{\infty}=0.$
\end{itemize}
\end{corollary}
Corollary \ref{equiv_compact} also follows from \cite[Proposition 2.3]{bulancea25}, where the equivalence between recurrence and uniform rigidity is obtained for multipliers in commutative unital Banach algebras. In \cite[Theorem 7.1]{costamanopari12} the result is obtained for the case when $K$ is connected and metric.

If $C_0(\Omega)$ is separable then $\Omega$ must be $\sigma$-compact and  a locally compact $\sigma$-compact space $\Omega$ is \emph{hemicompact}, that is, it contains  a fundamental family of compact sets $(K_n)_n$ such that for any compact subset $K\subseteq \Omega$ there is $n$ such that $K\subseteq K_n$. In this case   recurrence and rigidity are equivalent concepts:

 \begin{corollary}\label{rec=rigc0}
 	Let $\Omega$ be $\sigma$-compact and $\phi\in C_b(\Omega)$. Then $M_\phi:C_0(\Omega)\to C_0(\Omega)$ is recurrent if and only if it is rigid.
 \end{corollary}	

 \begin{proof}
 	Let  $(K_n)$ be an increasing fundamental sequence of compact subsets of $X$. If $M_\phi$ is recurrent, then $\phi(K_n)$ is a  Dirichlet subset of $\T$  for each $n\in\N$. Hence we can choose an increasing sequence $(m_n)_n$ such that  $\max_{z\in K_n}|\phi^{m_n}(z)-1|<1/n$ for each $n\in\N$. Therefore, $\lim_n \sup_{z\in K_j}|\phi^{m_n}(z)-1|=0$ for each $j\in\N$. We conclude by Proposition \ref{nonunitalcarac} (ii).
 \end{proof}

When $\Omega$ is discrete,  (i) and (ii)  in Proposition \ref{nonunitalcarac} and Corollary \ref{rec=rigc0} reduce to  the following result. The countable case can already be found in  \cite[Theorem 5.4]{costamanopari14}.
\begin{corollary}\label{numerable}	
	Let $\Omega$ be a   discrete space and  $\phi\in  C_b(\Omega)$.
\begin{itemize}
	\item[(i)]  $M_{\phi}:C_0(\Omega)\rightarrow C_0(\Omega)$ is recurrent if and only if $\phi(\Omega)\subseteq\T$.
	\item[(ii)]  $M_{\phi}:C_0(\Omega)\rightarrow C_0(\Omega)$ is rigid if and only if $\phi(\Omega)$ is Arbault.
	\item[(iii)] If $\Omega$ is countable, then $M_{\phi}:C_0(\Omega)\rightarrow C_0(\Omega)$ is recurrent if and only if it is rigid.
		\item[(iv)] If $\Omega$ is finite, then $M_{\phi}:C_0(\Omega)\rightarrow C_0(\Omega)$ is recurrent if and only if it is uniformly rigid.
\end{itemize}
\end{corollary}

Corollary \ref{numerable} implies that if  $\Omega$ is discrete and $\phi(\Omega)=\T$, then $M_{\phi}$ is recurrent but it is not rigid, as $\T$ is not an Arbault set. Therefore, the next example gives a negative answer for the analogue of \cite[Question 9.5]{costamanopari14} in the context of non-separable Banach spaces.
%

\begin{example}
	Let $I$ be an index set of continuum cardinal. Let $\T=\{\lambda_i:\ i\in I\}$ and $\phi:I\to \T, i\mapsto \lambda_i$. The operator $M_\phi:\ c_0(I)\to c_0(I)$ is recurrent but not rigid.
\end{example}


In order to get an example $M_{\phi}:C_0(\Omega)\rightarrow C_0(\Omega)$ that is rigid but not uniformly rigid,  it is enough to consider $\Omega$  countable and $\phi(\Omega)$ dense in a subset of $\T$ with positive measure. Since  the closure of a Dirichlet set is  a Dirichlet set, and a subset of $\T$ with positive measure cannot be Dirichlet, Corollary \ref{numerable} applies. We can, for instance, give the following example, which was already recorded in \cite[Example 5.5]{costamanopari14}.

\begin{example}
Let $\lambda=e^{2\pi\theta i}\in \T$, $\theta\in \R\setminus \Q,$ and $\phi:\Z\to \T, n\mapsto \lambda^n$. The operator $M_\phi:\ c_0(\Z)\to c_0(\Z)$ is rigid but not uniformly rigid.
\end{example}

%
%
%

 \begin{proposition}\label{phiconstant_C0}
 	Let $\Omega$ be  a locally compact Hausdorff space  	 and $\phi\in C_b(\Omega)$. Consider the following properties:
 	\begin{itemize}
 		\item[(i)] $M_{\phi}:C_0(\Omega)\rightarrow C_0(\Omega)$ is recurrent.
 			\item[(ii)]  $\phi(\Omega)\subseteq \T$ and $\phi$ is constant in each connected component of $\Omega.$
 	\end{itemize}
 If  $\Omega$ is locally connected or  the connected components of $\Omega$ are $\sigma$-compact , (i) implies (ii). If the connected components of $\Omega$ are  open, then (ii) implies (i). As a consequence, if $\Omega$ is locally connected or metrizable with open connected components, then (i) and (ii) are equivalent.
 \end{proposition}

 \begin{proof}
 Assume first that  (i) holds. First, $\phi(\Omega)\subseteq \T$ by Proposition \ref{lem:T_C0}.

  If $\Omega $ is locally connected and $C$ is a connected component of $\Omega$,  then $\phi^{-1}(\{a\})$ will be open for every $a\in \phi(C)$. Indeed, for every $x\in \phi^{-1}(\{a\})$ and  any connected relatively compact  neighbourhood $V$ of $x$, $\phi(V)$ will be, by Proposition \ref{nonunitalcarac}, a connected  Dirichlet set, i.e. $V\subseteq \phi^{-1}(\{a\})$.
Being open, and closed, $\phi^{-1}(\{a\})$ must contain $C$.

 If the connected components of $\Omega$ are $\sigma$-compact and  $\phi$ is non-constant on one of them, say $C$,  $\phi(C)$ is an interval in $\mathbb{T}$ and then  it has  the Baire property.  For a covering  $C=\bigcup_n K_n  $  by countably many compact sets, the interior of $\phi(K_n)$ will have to     be non-empty for some $n$. This contradicts Proposition \ref{nonunitalcarac}, since Dirichlet sets cannot contain non-empty open sets.

  %
%
   Assume now  that the connected components of $\Omega$ are  open and that  $\phi(\Omega)\subseteq \T$ and is constant in each connected component of $\Omega.$ Let us see that the  elements in  $C_0(\Omega)$ with  compact support, which form a dense subset in $C_0(\Omega),$  are recurrent vectors. Let $f\in C_0(\Omega)$, $\text{supp}(f)\subseteq K$, with $K\subseteq \Omega$ compact. The connected components of $\Omega$  are open, so there exist  $(U_j)_{j=1}^N$  connected components   such that $K\subseteq \bigcup_{j=1}^NU_j$.   Let $\alpha_j=\phi(U_j)\in \T,$ $1\leq j\leq N$. Thus, for fixed $n\in \N,$
 	$$\|M_\phi^n f-f\|_{\infty}=\|(\phi^n-1)  \sum_{j=1}^{N}  \Cf_{U_j} f\|_{\infty}=\|\sum_{j=1}^{N} (\alpha_j^n-1)\cdot  \Cf_{U_j} f\|_{\infty}\leq \|f\|_{\infty}\max_{1\leq j\leq N}|\alpha_j^n   -1|.$$ 	
 	Since the finite set  $(\alpha_j)_{j=1}^N$ is Dirichlet,  we obtain that  $f$ is a recurrent vector.

We finally observe that, in the setting of locally compact spaces, connected and metrizable spaces are $\sigma$-compact.
 \end{proof}


Regarding rigidity, Propositions \ref{rangofinito}, \ref{phiconstant_C0} and \ref{nonunitalcarac} yield the following:

\begin{proposition}
	\label{cc}
	Let $\Omega$ be a  locally compact Hausdorff space and let $\phi\in C_b(\Omega)$. Assume that $\Omega$ is  locally connected or  has open and $\sigma$-compact  connected components.
	\begin{itemize}
				\item[(i)]  $M_{\phi}:C_0(\Omega)\rightarrow C_0(\Omega)$ is rigid if and only if $\phi(\Omega)\subseteq \T$  is Arbault and it is constant in each connected component of $\Omega.$
		\item[(ii)] If $\Omega$ has a  countable number of connected components, then $M_{\phi}:C_0(\Omega)\rightarrow C_0(\Omega)$ is recurrent if and only if it is rigid.
		\item[(iii)] If $\Omega$ has a  finite number of connected components, then $M_{\phi}:C_0(\Omega)\rightarrow C_0(\Omega)$ is recurrent if and only if it is uniformly rigid if and only if $\phi(\Omega)\subseteq \T$ is finite.
	\end{itemize}
\end{proposition}
Proposition \ref{cc}(iii) yields the following extension of  \cite[Theorem 7.1]{costamanopari14}.
\begin{corollary}
	\label{finitock}
		Let $K$ be a   compact Hausdorff space with open connected components and let $\phi\in C(K)$. The multiplier $M_\phi:C(K)\to C(K)$ is recurrent if and only if it is uniformly rigid if and only if $\phi(K)\subseteq \T$ is finite.
\end{corollary}

 We see   below that when the connected components of a compact set $K$ are not open, then the range of a recurrent (uniformly rigid) multiplier does not need to be finite.

 \begin{example}
 \label{Cantor}
{\em Let $K\subseteq [0,1]$ be the Cantor set. There exists $\phi:  K\to\T$ continuous with infinite range such that $M_{\phi}:\ C(K)\to C(K)$ is uniformly rigid. We consider a Dirichlet subset of $\T$, say  $S=\{z_n=e^{\frac{2\pi i}{n!}}\colon \ n\in\N\}\subset \T$. Define $\phi$ by  $\phi(x)=z_n$ when  $x\in [1-1/3^{n-1},1-1/3^{n-1}+1/3^n]=[1-1/3^{n-1},1-2/3^{n}]$, $\phi(1)=1$.   Then $\lim_k\|M_{\phi}^{k!}-\mathbf{1}\|=0$ by Corollary \ref{equiv_compact}, since $S$ is Dirichlet.
}
\end{example}

%
%

\section{Multipliers on semisimple commutative Banach algebras}\label{Section:BanachAlgebras}

Having characterized recurrence and rigidity for multipliers of $C_0(\Omega)$, we aim to extend this analysis to multipliers of general semisimple commutative Banach algebras. 

Throughout this section, we assume that $\mA$ is a semisimple commutative Banach algebra. While isolated results hold more generally, this blanket assumption is typically required for our statements to be meaningful or valid.
	


We will find out that  the behaviour of  $T:\mA\rightarrow \mA$ is  strongly connected to the behaviour of the associated multiplication operator $\Mg{T}: C_0(\Delta(\mA))\rightarrow C_0(\Delta(\mA))$, even equivalent under mild assumptions.

 Proposition \ref{lem:T} shows that the condition $\widehat{T}(\Delta(\mA))\subseteq \T$ is necessary for $T$ to be recurrent;  this fact will be used throughout without further mention.

\begin{proposition}
	\label{densrecurrent}
Let $\mA$ be a commutative Banach algebra and let $T\in M(\mA)$.
If $T$ is uniformly rigid, then so is $M_{\widehat{T}}$ on
$C_0(\Delta(\mA))$. Moreover, if $\Gamma(\mA)$ is dense in
$C_0(\Delta(\mA))$, then recurrence and rigidity of $T$ are
inherited by $M_{\widehat{T}}$.
\end{proposition}
\begin{proof}
The assertion of uniform rigidity   is a direct consequence of $\|\widehat{T}^n-1\|_\infty \leq \|T^n-I\|$, which follows from the continuity of $\Gamma_M$.	
	
It is immediate from $\|\widehat{T}^n\widehat{a}-\widehat{a}\|_\infty \leq \|T^n(a)-a\|$ for every $a\in \mA$, that $\widehat{a}$ is a recurrent vector  for $M_{\widehat{T}}$ whenever $a$ is a recurrent vector for $T$. Hence the assertion concerning  recurrence follows immediately from the density of $\Gamma(\mA)$ in $C_0(\Delta(\mA))$. The same argument proves that $(M_{\widehat{T}^{n_k}})_k$ is a bounded sequence in $\mathcal{B}(C_0(\Delta(\mA)))$ which is convergent to $I$ in $\Gamma(\mA)$  whenever  $T$ is rigid along $(n_k)_k$.
The claim then  follows  from the density of $\Gamma(\mA)$ in $C_0(\Delta(\mA))$.
\end{proof}

\begin{remark}
	\label{discreto}
	By the Stone--Weierstrass theorem, the density condition in Proposition~\ref{densrecurrent} is satisfied whenever $\mA$ is self-adjoint. If $\Delta(\mA)$ is discrete, this condition holds even  without requiring $\mA$ to be semisimple: \v{S}ilov's Idempotent Theorem ensures that $\Cf_{\chi}\in\Gamma(\mA)$ for every $\chi\in\Delta(\mA)$.
\end{remark}

The following result, due to T. Körner, shows that the Dirichlet property entails approximation not only in the supremum norm but also in the algebra norm. For each $n\in \Z$, $\chi_n$ will denote the character $\chi_n \colon \T\to \T$ given by $\chi_n(t)=t^n,$ $t\in \T$.

Recall that the Fourier algebra $A(\T)$, also known in this case as the Wiener algebra,  consists of all functions $f=\sum_{n\in \Z} a_n \chi_n$ on $\T$ with $(a_n)_{n\in\Z} \in \ell^1(\Z)$, equipped with the norm $\Vert{}f\Vert{}_{A(\T)} = \sum_{n\in \Z}\vert{}a_n\vert{} < \infty$.

\begin{lemma}[Lemma 4.1 of \cite{korn73i}]\label{korner}
	Let  $E$ be a Dirichlet  subset of $\T$. If   $(n_k)_k$ satisfies  $\sup_{t\in E}\{|t^{n_k}-1|\}\rightarrow 0$ as $k\rightarrow \infty,$ then there exists a sequence $(\psi_{n_k})_k\subseteq A(\T)$ such that $\restr{\psi_{n_k}}{E}=\chi_{n_k}-1$ and $\|\psi_{n_k}\|_{A(\T)}\rightarrow 0$ as $k\rightarrow \infty.$
\end{lemma}
We first see how to exploit Lemma  \ref{korner}.
\begin{lemma}\label{recvec}%
Assume $\mA$ is self-adjoint. Let $ T\colon \mA\to \mA$ be  a power-bounded multiplier such that $\widehat{T}(\Delta(\mA))\subseteq \T$.  Let   $C\subseteq \Delta(\mA)$  be such that $\widehat{T}(C)$ is a Dirichlet subset of $\T$ along the sequence $(n_k)_k$.
Then
   $$\lim_k  \sup\{\|T^{n_k}(a)-a\|: \, a\in \mA, \  \supp{\widehat{a}}\subseteq C,\ \|a\|\leq 1\}=0.$$
In particular, every $a\in\mA$ satisfying
$\supp{\widehat{a}}\subseteq C$ is recurrent along $(n_k)_k$.
	
\end{lemma}
\begin{proof}
Let $\psi\in A(\T)$  and  $T\in M(\mA)$  satisfying the hypotheses. Since $\mA$ is self-adjoint, Proposition \ref{Remark_inversa} yields that $T^{-1}$ exists and it is also power bounded. Let $L:=\sup_{n\in\Z}\|T^n\|$.  We observe first that  the composition operator $R_T: A(\T)\mapsto M(\mA)$, $\psi(t) = \sum_{n \in \mathbb{Z}} a_n^\psi t^n \mapsto \sum_{n \in \mathbb{Z}} a_n^\psi T^n$ is well defined and bounded, with

\begin{equation}
	\label{keybound}
	\|R_T(\psi)\|\leq L \|\psi\|
	\end{equation}

	 From the definition we get that $R_T(\psi)=\sum_{n \in \mathbb{Z}} a_n^\psi T^n$ is a convergent series in $M(\mA)$.  	Continuity of $\Gamma_M$ yields

\begin{equation}
	\label{semisimple}
\widehat{R_T(\psi)}=	\Gamma_M(R_T(\psi))=\sum_{n \in \mathbb{Z}} a_n^\psi \widehat{T}^n=\psi\circ \widehat{T}.
\end{equation}

  By hypothesis, there exists a sequence $(n_k)_k$ such that $\widehat{T}^{n_k}-1\rightarrow 0$ uniformly on $C$. Lemma \ref{korner},  produces then a  sequence $(\psi_{n_k})_k\subseteq  A(\T)$  such that $\restr{\psi_{n_k}}{\widehat{T}(C)}=\restr{\chi_{n_k}-1}{\widehat{T}(C)}$ and $\|\psi_{n_k}\|_{\mbox{\tiny $A(\T)$}}\rightarrow 0$ as $k\rightarrow \infty.$  Fix $a\in \mA$ such that $\supp{\widehat{a}}\subseteq C.$

      We use \eqref{semisimple} to get
\begin{eqnarray*}
	\widehat{T^{n_k}(a)-a}=\widehat{T}^{n_k}\widehat{a}-\widehat{a}=((\chi_{n_k}-1)\circ\widehat{T})  \widehat{a}=(\psi_{n_k}\circ \widehat{T})\widehat{a}
	=\widehat{R_T(\psi_{n_k})}\widehat{a}=\widehat{R_T(\psi_{n_k})(a)}.
\end{eqnarray*}
Since $\mA$ is semisimple, $T^{n_k}(a)-a=R_T(\psi_{n_k})(a),$ thus we get from \eqref{keybound}
\begin{eqnarray*}
	\|T^{n_k}(a)-a\|	=\|R_T(\psi_{n_k})(a))\| \leq   L \|\psi_{n_k}\|_{\mbox{\tiny $A(\T)$}}\, \|a\|.
\end{eqnarray*}
This proves the Lemma,  since $\|\psi_{n_k}\|_{_{A(\T)}}\rightarrow 0$ as $k\rightarrow \infty$.
    \end{proof}




 \subsection{Recurrence}
Throughout  this section, we assume that $\mA$ is self-adjoint.  This condition   can be replaced in every case by requiring the multiplier $T$ to be bi-power-bounded,
that is, $\{\|T^n\|:n\in\mathbb{Z}\}$ is bounded, and  replacing
 the multiplication operator
$\Mg{T}\colon C_0(\Delta(\mA))\to C_0(\Delta(\mA))$ with  the  multiplication operator
$\Mg{T}\colon \Gamma_M(\mA)\to\Gamma_M(\mA)$.
The assumption of bi-power-boundedness is natural in this context:
by \cite[Proposition 3.2(ii)]{costamanopari14}, every recurrent
power-bounded operator is invertible and its inverse is
power-bounded.

\begin{theorem}\label{recA}
		Assume that $\mA$ is self-adjoint     and let $T: \mA\rightarrow \mA$ be a multiplier. Consider the following statements:
	\begin{enumerate}
		\item[(i)] $T: \mA\rightarrow \mA$  is recurrent.
		\item[(ii)]   $\Mg{T}:C_0(\Delta(\mA))\rightarrow C_0(\Delta(\mA))$ is  recurrent.
		\item[(iii)]   $\widehat{T}(\Delta(\mA))\subseteq \T$ and  $\widehat{T}(K)$ is Dirichlet for every compact $K\subseteq \Delta(\mA).$
	\end{enumerate}
Then, 	(i) $\Rightarrow$ (ii) $\Leftrightarrow$ (iii) are always satisfied. If $\mA$ is assumed to be Tauberian and $T$ is power bounded,    all the statements are equivalent.
\end{theorem}

\begin{proof}
(i)$\Rightarrow$(ii) follows by Proposition \ref{densrecurrent}.
\\
(ii) $\Leftrightarrow$ (iii) follows from  Proposition \ref{lem:T_C0} and  Proposition \ref{nonunitalcarac}(i).\\
To prove that 	(iii) $\Rightarrow$ (i), let $a\in \mA$ be such that $\widehat{a}$ has compact support $K\subseteq \Delta(\mA)$. Then, $\widehat{T}(K)$  is a Dirichlet set and  Lemma  \ref{recvec} proves that $a$ is a recurrent vector. Since the algebra $\mA$ is Tauberian (and so, elements whose  Gelfand transform has  compact support are dense in $\mA$)  we obtain that $T$ is recurrent.
\end{proof}

%


If $\Delta(\mA)$ is  locally connected, or has open and $\sigma$-compact connected components, power boundedness is not required for the equivalence of the three properties in Theorem \ref{recA}. Power boundedness, however, remains an essential assumption if the connected components of $\Delta(\mA)$ fail to be open (see Theorem \ref{C0_open_equiv_rec} in Subsection \ref{subs:abelian}). The Tauberian character of $\mA$  cannot be omitted either, even if $\Delta(\mA)$ is discrete and $T$ is power bounded, see Example \ref{cw}.

\begin{theorem}\label{rec_equiv_loc}
Assume that  $\mA$ is  Tauberian, $\Gamma(\mA)$ is dense in $C_0(\Delta(\mA))$, and     $\Delta(\mA)$  is either locally connected or  has open and $\sigma$-compact connected components. Let $T: \mA\to \mA$ be a multiplier.	
	The following statements are equivalent:	
	\begin{enumerate}
		\item[(i)] $T: \mA\rightarrow \mA$   is recurrent.
		\item[(ii)] 	$\Mg{T}:C_0(\Delta(\mA))\rightarrow C_0(\Delta(\mA))$ is recurrent.
		\item[(iii)]  $\widehat{T}(\Delta(\mA))\subseteq \T$ and it is constant in each connected component.
	\end{enumerate}
\end{theorem}

\begin{proof}
	
	(i)$\Rightarrow$(ii) follows from Proposition \ref{densrecurrent}. (ii)$\Leftrightarrow$(iii) from Proposition \ref{phiconstant_C0}.

Let us see that (iii)$\Rightarrow$(i)  by proving that elements in  $\mA$ whose Gelfand transform has compact support  are recurrent vectors for $T.$ Let $a\in \mA$ with $\text{supp}(\widehat{a})\subseteq K$, $K\subseteq \Delta(\mA)$ compact.
Under either hypothesis, the connected components of $\Delta(\mA)$ are open. Hence there exist  $(U_j)_{j=1}^N$   connected components such that $K\subseteq \bigcup_{j=1}^NU_j$.  We have  that  for $1\leq j\leq N$ there is $\alpha_j\in\C$ such that $\widehat{T}|_{U_j}=\alpha_j$. So, the span  of $\{\widehat{T}^n\widehat{a}: n\in\N_0\}$ is contained in that of $\{\Cf_{U_j}\widehat{a}: \ 1\leq j\leq N \}$, a finite dimensional subspace of $C_0(\Delta(\mA))$. Injectivity of the Gelfand transform yields that  $\{T^n(a):\  n\in\N_0\}$ spans also a finite dimensional subspace in $\mA$.   Therefore there exists $M>0$ such that
$$\|T^{n}(a)-a\|\leq M\|\widehat{T}^n\widehat{a}-\widehat{a}\|_\infty\leq M \sum_{j=1}^{N} |\alpha_j^n-1| \|\widehat{a}\|_\infty,$$

\noindent for all $n\in\N$ and then we conclude by the Dirichlet theorem.


\end{proof}

\begin{remark}
	\label{tauberianquasirigid}
	The proof of (iii)$\Rightarrow$(i) in Theorem~\ref{rec_equiv_loc}
	yields the following statement. Suppose that $\mA$ is Tauberian and
	that the connected components of $\Delta(\mA)$ are open. Let
	$T\colon\mA\to\mA$ be a multiplier such that
	$\widehat{T}(\Delta(\mA))\subseteq\T$ is a Dirichlet set. Since
	$\widehat{T}$ must be constant on each connected component of
	$\Delta(\mA)$, if  $(n_k)_k$  is the sequence satisfying
$	\lim_{k\to\infty}\sup\bigl\{|z^{n_k}-1|:
	z\in\widehat{T}(\Delta(\mA))\bigr\}=0,$
we get  $T^{n_k}(a)\to a$ for every $a\in\mA$ such that $\supp{\widehat{a}}$ is compact.  Hence, $T$ is quasirigid in the sense of \cite{GLP1,GLP2}.
\end{remark}

Taking into account Remark \ref{discreto}, the corresponding statement when the spectrum is discrete can be stated in a  cleaner way.

\begin{corollary}\label{Adiscr}
Assume that $\mA$ is Tauberian and that $\Delta(\mA)$ is  discrete. Let $T: \mA\rightarrow \mA$ be a multiplier.	
	The following statements are equivalent:	
	\begin{enumerate}
		\item[(i)] $T: \mA\rightarrow \mA$   is recurrent.
	\item[(ii)]	$\Mg{T}:c_0(\Delta(\mA))\rightarrow c_0(\Delta(\mA))$ is recurrent.
	\item[(iii)]$\widehat{T}(\Delta(\mA))\subseteq \T$.
	\end{enumerate}

\end{corollary}

Now we see that $\mA$ needs to be Tauberian.

\begin{example}\label{cw}
  There is a semisimple commutative  Banach algebra $\mA$  with $\Delta(\mA)$ discrete that admits a  power-bounded multiplier $T$ with $\widehat{T}(\Delta(\mA))\subseteq \T$ but $T$ is not recurrent.
\end{example}
\begin{proof}
We consider  here the algebra introduced by Mirkil \cite{mirk59-60} as a standard example of a non-Tauberian semisimple Banach algebra.
Let $\mathfrak{c}$ denote the Banach algebra of convergent sequences  with the supremum norm. We define  $\mA$ as     the semisimple Banach algebra obtained by equipping $\mathfrak{c}$ with the multiplication \[\mathbf{a}\cdot \mathbf{b}=\left(\alpha_n a_n b_n\right)_n,\]
for any sequences $\mathbf{a}=(a_n)_n,\; \mathbf{b}=(b_n)_n\in \mathfrak{c}$, where $\alpha_n=1/n$ for all $n\in\N$.  It is easy to show that $\Delta(\mA)=\{\alpha_n\delta_n: \ n\in\mathbb{N}\}$, which is a discrete space.
  Define an isometric multiplier $T$ of $\mA$ as $T(\mathbf{a})=(e^{ \frac{\pi i}{n}} a_n)_n $, for each $\mathbf{a}=(a_n)_n\in \mathfrak{c}$. The range of $\widehat{T}$ is contained in $\T$, since, for each $n\in \N$,
  \[\widehat{T}(\alpha_n \delta_n)=e^{ \frac{\pi i}{n}} .\]
 The multiplier $T$ has non-recurrent vectors: if  $\mathbf{1}$ is the constant 1-sequence  and $k\in \N$,
 \[\norm{ T^k(\mathbf{1})-\mathbf{1}}=\norm{\left(e^{ \frac{\pi k i}{n}}-1 \right)_n}=2\]
  Since $T$ is power bounded, it can only be recurrent when all vectors are recurrent.
\end{proof}

\subsection{Rigidity and uniform rigidity}\label{Section:BanachAlgebras_rig}

In the following result we characterize rigid power bounded multipliers on  commutative Banach algebras.

\begin{theorem}\label{rigidA}
	Let $\mA$ be Tauberian and self-adjoint, and let $T: \mA\rightarrow \mA$ be a power bounded multiplier. The following statements are equivalent:
	\begin{enumerate}
		\item[(i)] $T: \mA\rightarrow \mA$  is  rigid.
		\item[(ii)]   $\Mg{T}:C_0(\Delta(\mA))\rightarrow C_0(\Delta(\mA))$ is    rigid.
		\item[(iii)] There is an increasing sequence $(n_k)$ of natural numbers such that \[\lim_k\sup_{\chi\in K}|\widehat{T}^{n_k}(\chi)-1|=0\]
for every compact set $K\subseteq\Delta(\mA)$.
	\end{enumerate}
\end{theorem}
 \begin{proof}
 		(i)$\Rightarrow$(ii)	follows from Proposition \ref{densrecurrent}	\\
 		(ii)$\Leftrightarrow$(iii) is a consequence of Proposition  \ref{nonunitalcarac} (ii). \\
 		For (iii)$\Rightarrow$(i) assume that (iii) holds. Let $(n_k)$ be a sequence such that $\widehat{T^{n_k}-1}$ converges to 0 uniformly on every compact subset $K\subseteq \Delta(\mA)$. We apply Lemma  \ref{recvec} to obtain that $\lim_k \|T^{n_k}(a)-a\|=0$ for any $a$ such that $\widehat{a}$ is compactly supported. Since $\mA $ is Tauberian and $T$ is power bounded, we conclude  that $\lim_k \|T^{n_k}(a)-a\|=0$ for any $a\in \mA$.	
  \end{proof}


The next result follows directly from Corollary~\ref{rec=rigc0} and Theorem~\ref{rigidA}. Its hypotheses are satisfied whenever $\mA$ is separable, thereby answering \cite[Question 9.5]{costamanopari14} in the affirmative for multipliers on self-adjoint, Tauberian, semisimple commutative Banach algebras, as noted in the introduction.

 \begin{corollary} \label{sigmacompact}
 	Let $\mA$ be Tauberian and self-adjoint and let $T: \mA\rightarrow \mA$ be a power bounded multiplier.
Assume that   $\Delta(\mA)$ is $\sigma$-compact. A  power bounded multiplier of $\mA$ is recurrent if and only if it is rigid.
 \end{corollary}	

Unlike recurrence (cf. Corollary~\ref{Adiscr}), rigidity still requires in the discrete case the additional assumption of power boundedness to guarantee that  $T:\mA\rightarrow \mA$ and $\Mg{T}:C_0(\Delta(\mA))\rightarrow C_0(\Delta(\mA))$ behave equivalently.  Proposition \ref{abelian-block-construction} provides a counterexample showing that  this assumption cannot be omitted.

 \begin{corollary}\label{magdiscpbrig}
Assume that  $\mA$ is  Tauberian and self-adjoint, and  that  $\Delta(\mA)$ is discrete and let $T:\mA\rightarrow \mA$  be  a power bounded multiplier. The following statements are equivalent.
	\begin{enumerate}
		\item[(i)] $T: \mA\rightarrow \mA$  is  rigid.
		\item[(ii)] $\Mg{T}:C_0(\Delta(\mA))\rightarrow C_0(\Delta(\mA))$ is rigid.
		\item[(iii)] $\widehat{T}(\Delta(\mA))$ is an Arbault subset of $\T $.
	\end{enumerate}
	If $\Delta(\mA)$ is countable, then (i), (ii) and (iii) are equivalent to:
	\begin{enumerate}
		\item[(iv)] $\widehat{T}(\Delta(\mA))\subseteq \T$.
		\item[(v)] $T$ is recurrent.
		\item[(vi)]  $\Mg{T}:C_0(\Delta(\mA))\rightarrow C_0(\Delta(\mA))$ is recurrent.
	\end{enumerate}
\end{corollary}


For uniform rigidity, we get the next result:

\begin{theorem}\label{URA}
 Assume that  $\mA$ is  Tauberian. Let $T: \mA\rightarrow \mA$ be a multiplier. Consider the following statements:
 \begin{enumerate}
 	\item[(i)] $T: \mA\rightarrow \mA$  is uniformly rigid.
 	\item[(ii)]  $\Mg{T}:C_0(\Delta(\mA))\rightarrow C_0(\Delta(\mA))$ is   uniformly rigid.
 	\item[(iii)] $\widehat{T}(\Delta(\mA))\subseteq \T$ is Dirichlet.
 \end{enumerate}
 (i) $\Rightarrow$ (ii) $\Leftrightarrow$ (iii) is always satisfied.  If $\mA$ is self-adjoint and  $T$ is power bounded, then the three statements are equivalent.
\end{theorem}

\begin{proof}
(i) $\Rightarrow$ (ii) follows from Proposition \ref{densrecurrent}.\\
(ii) $\Leftrightarrow$ (iii) is Proposition \ref{nonunitalcarac}(iii).\\
	(iii) $\Rightarrow$ (i).
 It follows directly from Lemma \ref{recvec}.
\end{proof}

\begin{remark}
	\label{normality}
By Proposition~\ref{nonunitalcarac}, the density of $\Gamma(\mA)$ in $C_0(\Delta(\mA))$ in Proposition~\ref{densrecurrent} can be replaced by the following \emph{normality} condition: ``for each $K\subseteq \Delta(\mA)$, there exists $u\in \mA$ such that $\widehat{u}\vert{}_K=1$.'' This condition is satisfied by every regular Banach algebra. Theorem~\ref{rec_equiv_loc} can similarly be restated using this normality condition in place of density. Furthermore, as already emphasized at the beginning of Subsection 4.1, self-adjointness can be replaced by normality in Theorems~\ref{recA}, \ref{rigidA}, and \ref{URA} (as well as their corollaries), provided that $T$ is assumed to be bi-power-bounded (i.e., both $T$ and $T^{-1}$ are power-bounded).\end{remark}

%

\section{Convolution operators and the Fourier algebra}\label{Section:A(G)}

In this section, we consider the multipliers of the Fourier algebra $A(G)$
of a locally compact group $G$. This encompasses the multipliers of the
group algebra $L^1(G)$ of a locally compact Abelian group $G$ as a particular case.

For an arbitrary locally compact group $G$, the Fourier algebra $A(G)$ and the Fourier--Stieltjes algebra $B(G)$ are   subalgebras of the algebra of bounded continuous functions on $G$. The algebra $B(G)$ is unital, and $A(G)$ is an ideal of $B(G)$ contained in $C_0(G)$. Compactly supported functions are dense in $A(G)$, when $G$ is discrete $\ell_2(G)\subseteq A(G)$.

  One can introduce in $B(G)$ a suitable norm, whose definition we do not require here, that makes $B(G)$ into a Banach algebra, and makes $A(G)$ closed in  $B(G)$. When $G$ is Abelian, the Fourier and Fourier--Stieltjes transforms establish respective isometric algebra isomorphisms from the convolution algebras $L^1(\widehat{G})$ and $M(\widehat{G})$  onto the function algebras $A(G)$ and $B(G)$, respectively. For a brief introduction to these algebras, see Section~2.2 of \cite{galijordarodr25}. Eymard's paper \cite{eyma64} is the  reference to the  foundational work; see also \cite{rudin90} for the commutative setting, and \cite{kanilaubook} for a comprehensive exposition.

  We  need to bring up here   the
rich operator-theoretic structure  that Fourier and Fourier-Stieltjes algebras possess. It  largely stems  from their roles as preduals of important von Neumann algebras.  Here we will just touch upon the essentials concerning the relationship between $A(G)$ and the group von Neumann algebra $VN(G)$  that will be required in the sequel.

Recall that  $VN(G)$ is the  von Neumann algebra generated by the left-translation operators
$\{\lambda_x:x\in G\}$, where
\[
(\lambda_x f)(y)=f(x^{-1}y), \qquad f\in L^2(G),\ x,y\in G,
\]
 that is, $VN(G)$ is the ultraweak closure of the linear span of
$\{\lambda_x\colon x\in G\}$ in $\mathcal{B}(L^2(G))$.

According to \cite[Theorem 2.4.3]{kanilaubook}, the Fourier algebra $A(G)$ consists of  all functions $f\ast\tilde{g}$, where $f,\, g\in L^2(G)$ and $\tilde{g}(x)=\overline{g(x^{-1})}$ for every $x\in G$.
Alternatively \cite[Proposition~2.3.3]{kanilaubook}, $A(G)$ can be characterized as the closed linear span of
\[
\{f \ast \tilde{g} \colon f, g \in C_c(G)\},
\]
where $C_c(G)$ denotes the space of continuous functions with compact support.

There is an isometric isomorphism between $VN(G)$ and the dual Banach space $A(G)^\ast$ of $A(G)$. Under this identification, each operator $S \in VN(G)$ corresponds to a functional $\phi_S \in A(G)^*$  so that, for any $u=f\ast \tilde{g},  f, g \in C_c(G)$,
\[   \<{\phi_S,u}=\bigl(S(\overline{g})\ast \check{f}\bigr)(e),\]
where $\check{f}(x)=f(x^{-1})$, see  \cite[Page 55]{kanilaubook}.

 It follows that (see \cite[Theorem 2.3.9 (ii)]{kanilaubook}$) \<{\lambda_x,u}=u(x)$ for every $x\in G$ and $u\in A(G)$.

The spectrum $\Delta(A(G))$ of $A(G)$ is naturally identified with $G$ \cite[Theorem 2.3.8]{kanilaubook}. Under this identification, the Gelfand transform on $A(G)$ becomes the identity. Every multiplier $T$ on $A(G)$  acts therefore as a multiplication operator on $C_0(G)$, $M_\phi(u) = \phi \cdot u$, where $\phi = \widehat{T}\in C_b(G)$. Conversely, any function $\phi \in C_b(G)$ satisfying $\phi \cdot A(G) \subseteq A(G)$ defines a multiplier $M_\phi$ on $A(G)$. With a slight abuse of notation, we denote by $M(A(G))$ the algebra of all such functions.
Accordingly,  if $\phi\in M(A(G))$, $M_\phi$ is a multiplier of $A(G)$ and the Gelfand transform of $M_\phi$  is precisely the function $\phi$.
 Since  $A(G)$ is an ideal in $B(G)$ we always have that $B(G)\subseteq M(A(G))$.


In Section \ref{Section:BanachAlgebras}, we compared the operators $T\colon \mA\to \mA$ and $\Mg{T}\colon C_0(\Delta(\mA)) \to C_0(\Delta(\mA))$. When $T=M_\phi$, both $T$ and $\Mg{T}$ are multiplication operators with symbol $\phi$. In this case, in the context of the Fourier algebra, instead of $\Mg{T}$, we will prefer the notation $M_\phi^0 \colon C_0(G)\to C_0(G)$ for this operator. Evidently, $\restr{M_\phi^0}{A(G)}=M_\phi$. The primary distinction between $M_\phi^0$ and $M_\phi$ lies in the norm of the underlying domain, that of $A(G)$ for $M_\phi$ and the supremum norm  for $M_\phi^0$. It must be noted that for every $u\in A(G)$, $\norm{u}_{C_0(G)}=\norm{u}_\infty\leq \norm{u}_{B(G)}$, and that  these norms fail to be equivalent, unless $G$ is finite. When $G$ is amenable (in particular, when $G$ is abelian or compact), $M(A(G))=B(G)$, in which case $\norm{\phi}_{B(G)}=\norm{M_\phi}$ for any $\phi\in B(G)$. If $G$ is nonamenable, however, $M(A(G))$ strictly contains $B(G)$, the inequality $\norm{M_\phi} \leq \norm{\phi}_{B(G)}$ holds, and the multiplier norm and the $B(G)$ norm are not equivalent, see \cite[Chapter 5]{kanilaubook}.

The Fourier algebra is a  commutative, semisimple, regular, self-adjoint, and Tauberian Banach algebra, see   \cite[Proposition 3.8, p.210]{eyma64} or  \cite[Theorem 2.10.2]{kanilaubook}. Consequently, the results on recurrence and rigidity established in Section~\ref{Section:BanachAlgebras} apply directly to the operators $M_\phi$ and $M_\phi^0$.

           In the case of the Fourier algebra, the dynamics of $M_\phi$, $\phi\in M(A(G))$ are governed by the intricate interplay between the multiplier $\phi$ and the algebraic and topological properties of the group. To begin with,  we know from  Proposition \ref{densrecurrent} that $M_\phi^0$ is recurrent, rigid, or uniformly rigid whenever $M_\phi$ is. Hence  we get that  $\phi( G)\subseteq  \T$ when $M_\phi$ is recurrent. This is  a demanding requirement, as will be seen.

      We have organized this section into three main parts. Firstly, we examine the algebraic properties that are forced  upon a power-bounded multiplier $\phi$ under the assumption that  $\phi(G)\subseteq \T$  (in which case we say that $\phi$ is \emph{unimodular}).  Secondly, we analyze the recurrence and (uniform) rigidity of these operators. Finally, we  see how the dynamics become less constrained when the  multipliers are not power bounded, and how they can even be used to break the equivalences established in the preceding sections.

\subsection{Algebraic properties of unimodular multipliers}

An element of $B(G)$ taking values in $\mathbb{T}$ must satisfy strong
algebraic constraints, particularly when its norm is small. We begin
by illustrating this phenomenon for elements of norm one.
%
%
%


\begin{lemma}\label{iso}
	Let     $\phi \in M(A(G))$ with $\phi(G)\subseteq \T$ and $\norm{M_\phi}\leq 1$. Then the adjoint operator $M_\phi^\ast \colon VN(G)\to VN(G)$ is a surjective isometry.
\end{lemma}
\begin{proof}
	It is enough to see that $M_\phi$ is a surjective isometry.  This is a simple consequence of the inverse $M_{\bar{\phi}}$ being again a multiplier of the same norm, which follows from the fact that for every $u\in A(G)$, we have $\bar{u}\in A(G)$ and $\norm{u}=\norm{\bar{u}}$ \cite[Remark 2.1.10(3)]{kanilaubook}).
\end{proof}

\begin{theorem}[Kadison, Theorem 7 of \cite{kadi54}]\label{kadi}
	Let $\mA$ and $\mathcal{B}$ be unital $C^\ast$-algebras and let $\Psi\colon \mA\to \mathcal{B}$ be a surjective linear isometry with $\Psi(I)=I$. Then, for every $S,T\in \mA$,
	\[\Psi(ST+TS)=\Psi(S)\Psi(T)+\Psi(T)\Psi(S).\]\end{theorem}
\begin{theorem}
	\label{homom}
	If      $\phi \in M(A(G))$ with $\phi(G)\subseteq \T$ and $\norm{M_\phi}\leq 1$, then $\phi$ is affine.
\end{theorem}
\begin{proof}
It suffices to prove that $\phi$ is a homomorphism under the additional assumption that $\phi(e)=1$, for, otherwise we work with $\tilde{\phi}=\phi\cdot \phi(e)^{-1}$.

  For $x\in G$, we have that	$\lambda_x\in VN(G)$, acts as  $ \<{\lambda_x,u}=u(x)$  \cite[Theorem 2.3.9 (ii)]{kanilaubook} for every $u\in A(G)$. Hence,  $M_\phi^\ast \lambda_x=\phi(x) \lambda_x$. In particular,  $M_\phi^\ast \lambda_e=\lambda_e.$  Theorem  \ref{kadi} then applies to the map $M_\phi^\ast\colon VN(G)\to VN(G)$. Applying this identity to the elements $\lambda_x, \lambda_{y^{-1}} \in VN(G)$ and noting that $\lambda_x \lambda_y = \lambda_{xy}$, we obtain
	\begin{align*}
		\phi(xy^{-1})\lambda_{xy^{-1}}+\phi(y^{-1}x)\lambda_{y^{-1}x}&=\phi(x)\cdot \phi(y^{-1})\left(\lambda_{x y^{-1}}+\lambda_{y^{-1}x}\right).
	\end{align*}
	If $xy^{-1}=y^{-1}x$, it follows right away that $\phi(x)\phi(y^{-1})=\phi(xy^{-1})$.
		If $xy^{-1}\neq y^{-1}x$, pick $u\in A(G)$ with $u(x y^{-1})=1$ and $u(y^{-1}x)=0$, then
	\begin{align*}
		\phi(xy^{-1})&=\<{\phi(xy^{-1})\lambda_{xy^{-1}}+\phi(y^{-1}x)\lambda_{y^{-1}x},u}\\
		&=\<{\phi(x)\cdot \phi(y^{-1})\left(\lambda_{x y^{-1}}+\lambda_{y^{-1}x}\right),u}\\
		&=\phi(x) \phi (y^{-1}).
	\end{align*}
\end{proof}


Nonaffine functions $\phi \in B(\mathbb{Z})$ that take its values in $T$ are easy to come by. The next example just singles one out.

\begin{example}\label{exnonaff}
Consider the characters of $\Z$: $\chi_1(n) = i^n$ and $\chi_2(n) = (-i)^n$. Define the element of $B(\Z)$,
	\[ \phi = \left( \frac{1-i}{2} \right) \chi_1 + \left( \frac{1+i}{2} \right) \chi_2. \]
Then $\phi(G)\subseteq \T$, $\norm{\phi}=\sqrt{2}$ and $\phi$ is not affine.
\end{example}
\begin{proof}
Direct computation yields
	\[ \phi(4k)=1, \quad \phi(4k+1)=1, \quad \phi(4k+2)=-1, \quad \phi(4k+3)=-1,\]
$\phi$ is clearly nonaffine, as $\phi(0)=1$ but $\phi(2)\neq \phi(1)\phi(1)$.

For the norm,  observe that $\phi=\widehat{\mu}$ with
 $\mu = \frac{1-i}{2} \delta_{-i} + \frac{1+i}{2} \delta_{i}$. Therefore,
\[ \|\phi\|_{\mbox{\tiny $B(\mathbb{Z})$}} = \|\mu\|_{M(\mathbb{T})} = \left| \frac{1-i}{2} \right| + \left| \frac{1+i}{2} \right|  = \sqrt{2}. \]	
\end{proof}

Theorem~\ref{homom} extends to every power bounded element of $B(G)$,
at the cost of weakening the algebraic properties of $\phi$; this is a consequence of  Lemma 3.1 of \cite{kanilauulger11}, which follows closely the Abelian proof of  \cite[Lemma 3.1]{schr70}.
To state the properties of the
resulting map, we need the notion of a piecewise affine map. We do not give its definition in the most natural form. Instead, we
incorporate some known characterizations that make it easier to use
(see \cite[Chapter~4]{rudin90} or \cite{schr70}; see also \cite{daws21}
for a thorough account of the subtleties surrounding piecewise affine
maps).

\begin{definition}\label{def:paf}
	Let $G$ and $H$ be  two topological  groups and let $\phi\colon G\to H$ be continuous. We say that $\phi$ is \emph{piecewise affine} if   $G=\bigcup_{j=1}^n S_j$ with
	\[ S_j=L_j \setminus \bigcup_{k=1}^{m_j} M_{j}^k,\]
	where $L_1,\ldots, L_n$ and $M_j^k$, $1\leq k\leq m_j$, $1\leq j \leq n$, are  open cosets in $G$ and there are affine maps
	$\alpha_j \colon L_j\to H$ such that $\restr{\phi}{S_j}=\restr{\alpha_j}{S_j}$.
\end{definition}
We can now  state the already mentioned extension of Theorem~\ref{homom}.

	\begin{theorem}[Lemma 3.1 of \cite{kanilauulger11}]	\label{paffdense+pbpa}
		If $\phi\in B(G)$ is power bounded and $\phi(G)\subseteq \T$, then $\phi$ is piecewise affine.
	\end{theorem}
%
%
%
%

\begin{lemma}\label{lem:coset}
	Let  $L$ be a coset in a group  $G$  and let $\alpha \colon L\to \T$ be an affine map.
	Consider $M_1,\ldots , M_n$ be cosets contained in $L$ with $L\neq \bigcup_{k=1}^n M_k$. If $\alpha (L)$ is dense in $\T$, then
	$\alpha\left(L\setminus \bigcup_{k=1}^n M_k\right)$ is also dense.\end{lemma}
\begin{proof}
In this proof we will use freely that subgroups (and, hence, cosets) of $\T$ are either finite or dense, and, hence, that all proper closed ones  have empty interior.
	
	Suppose first that the cosets  $\alpha(M_1),\ldots, \alpha(M_n)$ are all finite. Then
	\begin{align*}
		\overline{\alpha(L\setminus \bigcup_{k=1}^n M_k)}&\supset \overline{\alpha(L)\setminus \bigcup_{k=1}^n \alpha(M_k)}\\
		&=\overline{\alpha(L)}=\T.
	\end{align*}
	%
	
	We proceed by induction on $n$ assuming that $\alpha(M_1)$ is infinite.

Observe that for any $x\in M_k$, we can put  $L=xH$ and $M_k=xH_k$, for some subgroup $H$ of $G$ and some subgroup  $H_k$  of $H$.
	
	For  $n=1$, let $x\in M_1$. Take $y\in L\setminus M_1$ (which is not empty by hypothesis) and consider $M_1^\prime:=yH_1.$
We have that $M_1^\prime\subset yH=xH=L$ and $M_1^\prime\cap M_1=\emptyset$.
Being of the same cardinality as $\alpha(M_1)$,  we know that $\alpha(M_1')$ is dense  in $\T.$ Since $\alpha(M_1')\subset \alpha(L\setminus M_1)$, we see that $\alpha(L\setminus M_1)$ has to be dense.
	
	Assume now the result is true when the number of cosets involved is $n,$ and consider the case
	$L\setminus \bigcup_{k=1}^{n+1} M_k$. Now, take $y\in L\setminus \bigcup_{k=1}^{n+1} M_k$ (which is not empty by hypothesis) and define $M_1^\prime=yH_1$. As above, since $\alpha(M_1)$ is infinite, $\alpha(M_1^\prime)$ is dense.   Consider  $\alpha(M_1^\prime\setminus \bigcup_{k=2}^{n+1} (M_k\cap M_1'))$.
Since  $y\in L\setminus \bigcup_{k=1}^{n+1} M_k$, we have that  $M_1^\prime \neq \bigcup_{k=2}^{n+1} (M_k\cap M_1^\prime)$.
By the induction hypothesis, $ \alpha(M_1^\prime\setminus \bigcup_{k=2}^{n+1}( M_k\cap M_1^\prime ) $ is dense in $\T.$ Since
	\[ \alpha\left(M_1^\prime\setminus \bigcup_{k=2}^{n+1} (M_k\cap M_1^\prime)\right) \subseteq \alpha\left(L\setminus \bigcup_{k=1}^{n+1} M_j\right),\]
	the result follows.	
\end{proof}

\begin{corollary}
	\label{paffdense}
	Let $\phi\in B(G)$ be power bounded and $\phi(G)\subseteq \T$. Then $\phi(G)$ is either finite or dense in $\T$. Consequently, if $\phi(G)\subseteq \T$ is a Dirichlet set then $\phi(G)$ is finite. \end{corollary}

\begin{proof}
By Theorem \ref{paffdense+pbpa} we get that  $\phi$ can be described as in Definition \ref{def:paf}. If we assume that $\phi(G)$ is  infinite, then there exists $j\in \N$ with $1\leq j\leq n$ such that $\alpha_j(L_j)$ is infinite, and hence dense. Lemma \ref{lem:coset} shows that $\alpha_j\left(L_j\setminus \bigcup_{k=1}^{m_j} M_j^k\right)$ is dense in $\T$, hence the conclusion follows. If we assume $\phi(G)$ is Dirichlet, then also $\overline{\phi(G)}$ is so. Therefore $\phi(G)$ cannot be dense in $\T.$
\end{proof}

\begin{example}\label{exnonaffinf}
{\em   It is  easy to produce examples of power-bounded piecewise affine maps with values in $\T$, both with finite and infinite range.
	
The function $\phi \in B(\Z)$ introduced in Example \ref{exnonaff} has a finite image $\phi(\Z) \subseteq \T$, thus $M_{\phi}$ is power bounded and uniformly rigid in $A(\Z)$ by Proposition \ref{rangofinito}.


A slight modification produces a piecewise affine map with values in $\T$ with $M_\phi$ rigid but not uniformly rigid. 	Let $t_0 \in \mathbb{T}$ be an element of infinite order. Consider the characters $\chi_1, \chi_2 \in \widehat{\mathbb{Z}}$ defined by $\chi_1(n) = (i t_0)^n$ and $\chi_2(n) = (-i t_0)^n$, and define the function
\[ \phi = \left( \frac{1-i}{2} \right) \chi_1 + \left( \frac{1+i}{2} \right) \chi_2. \]
Direct evaluation at each $n \in \mathbb{Z}$ yields:
\[ \phi (4n)=t_0^{4n},\quad \phi(4n+1)= i t_0^{(4n+1)},\quad \phi(4n+2)=-t_0^{(4n+2)} \mbox{ and } \phi^{4n+3}=(-i) t_0^{(4n+3)}.\]
Since $t_0$ has infinite order, the image $\phi(\mathbb{Z})$ is an infinite dense subset of $\mathbb{T}$. 	 The measure associated with $\phi$ is $\mu = \frac{1-i}{2} \delta_{-it_0^{-1}} + \frac{1+i}{2} \delta_{it_0^{-1}}$, which satisfies $\|\phi\|_{B(\mathbb{Z})} = \sqrt{2} > 1$, as it should, as $\phi$ is not affine.
Since $\phi$ is a linear combination of two characters, a simple calculation using the binomial expansion shows that $\phi$ is power-bounded.}

\end{example}

The following result  shows that the  range of a function in $B(G)$ that is  power bounded in the multiplier norm can be an infinite Dirichlet subset of $\T$, in contrast with the case of power boundedness in the $B(G)$-norm, established in Corollary \ref{paffdense}. The construction relies on the existence of Leinert sets, and therefore  applies, for instance, to free groups on at least two generators.

By a characterization due to Bo\.zejko \cite{bozejko}, a subset
$E$ of a discrete group $G$ is a {\it Leinert set} if and only if
\[
\ell^\infty(E)\subseteq M(A(G)).
\]
\begin{proposition}
	\label{leinert}
Let $G$ be a discrete group which has an infinite Leinert  subset $E$. There exists $\phi\in B(G)$ such that $(\phi^k)_k$ is bounded in $ M(A(G))$ and $\phi(G)\subseteq \T$ is an infinite Dirichlet set.  	
\end{proposition}
\begin{proof}
Let us assume without loss of generality that $E=\{e_n:\ n\in\N\}$ is countable. Define $h_k:G\to\C$ in $M(A(G))$ as $h_k(e_n)=e^{(2\pi k i)/ n!}-1$ for each $n\in\N$, $h_k(g)=0$ for $g\notin E$. Define $\phi=1+h_1$. We get immediately $\phi^k=1+h_k$. Since $E$ is Leinert, the inclusion $\ell^\infty (E)\subseteq M(A(G))$  is continuous due to the closed graph theorem,  and then there is $C>0$ such that $\|f\|_{ M(A(G))}\leq C\|f\|_{\ell^\infty(E)}$ for each $f\in \ell^\infty(E)$. Hence we get $\|\phi^k\|_{ M(A(G))}\leq \|1\|_{ M(A(G))}+\|h_k\|_{ M(A(G))}\leq 1+C\|h_k\|_{\ell^{\infty}(E)}=1+2C$. Therefore $(\phi^k)_k$ is bounded in $M(A(G))$. The inclusion $\ell^2(G)\subseteq B(G)$
 and $\phi=h_1+1$ implies $\phi\in B(G)$.
\end{proof}
	
	\subsection{Recurrence and rigidity of power bounded  multipliers of $A(G)$.}
	As an immediate consequence of Theorems \ref{recA}, \ref{rigidA} and  \ref{URA}, the recurrence and rigidity properties of $M_\phi$ are equivalent to that of $M_\phi^0$ when $\phi\in M(A(G))$ is power bounded. Using Proposition \ref{nonunitalcarac},  Corollary \ref{equiv_compact} and Corollary \ref{rec=rigc0},  we can state this equivalence in the following way:
	
	\begin{theorem}
		\label{pb}
		Let $G$ be a locally compact group and let $\phi\in  M(A(G))$ be power bounded such that $\phi(G)\subseteq \T$.
		\begin{itemize}
			\item[(i)]  $M_{\phi}:A(G)\rightarrow A(G)$  is recurrent if and only if $\phi(K)$ is a Dirichlet set for every compact set $K\subseteq G$. When $G$ is $\sigma$-compact  $M_\phi$ is recurrent if and only if it is  rigid.
			\item[(ii)] $M_{\phi}:A(G)\rightarrow A(G)$  is  rigid if and only if there is an increasing sequence $(n_k)$ of natural numbers such that $\lim_k\sup_{x\in K} |\phi^{n_k}(x)-1|=0$  for every compact set $K\subseteq G$.
			\item[(iii)] 	$M_{\phi}:A(G)\rightarrow A(G)$  is uniformly  rigid if and only if $\phi(G)$ is Dirichlet.
			\item[(iv)] If $G$ is compact then $M_{\phi}:A(G)\rightarrow A(G)$  is recurrent if and only if $\phi(G)$ is Dirichlet.
		\end{itemize}
	\end{theorem}
	
	For uniformly rigid operators more can be said about the range of $\phi$ if we assume $\phi\in B(G)$ is power bounded. We remark that the proof of the next Theorem does not rely on Lemma \ref{korner}.
	
	\begin{theorem}
			\label{amenableur}
		Let $G$ be a locally compact group and let $\phi\in B(G)$ be power bounded with respect to $\|\ \|_{B(G)}$. Then $M_\phi$ is uniformly rigid if and only if $\phi(G)$ is finite.
	\end{theorem}
	
	\begin{proof}
		If $M_\phi$ is uniformly rigid then also $M_\phi^0$ is so and $\phi(G)\subseteq \T$ is Dirichlet by Theorem \ref{URA}. 	Corollary \ref{paffdense} shows that $\phi(G)$ is finite. The converse is also true by Proposition \ref{rangofinito}.
	\end{proof}	
	
	\begin{remark}
		Proposition \ref{leinert} and 	Theorem \ref{pb} show that Theorem~\ref{amenableur} needs not hold  for functions $\phi \in B(G)$ that are only  power-bounded in the multiplier norm.
	\end{remark}
	
	\subsection{Beyond Power Boundedness}
	
	\subsubsection{Discrete and compact case}
		For the Fourier algebra,  Theorem \ref{rec_equiv_loc}  adopts a more general form, as connected locally compact groups are always $\sigma$-compact.
	
	\begin{theorem}\label{rec_equiv_locA}
		Let $G$ be a locally compact group with open connected component of the identity and $\phi\in  M(A(G)).$
		The following statements are equivalent:
		\begin{enumerate}
			\item[(i)] $M_\phi$ is recurrent.
			\item[(ii)]  $M_\phi^0$ is recurrent.
			\item[(iii)]  $\phi(G)\subseteq \T$ and $\phi$  is constant in each connected component.
		\end{enumerate}
	\end{theorem}We will prove in Theorem \ref{C0_open_equiv_rec} below that the equivalence in Theorem \ref{rec_equiv_locA} fails in every group whose identity component is not open. This will also yield a multiplier $M_\phi$ (necessarily not power-bounded)  such that $M_\phi^0$ is recurrent (even uniformly rigid), but $M_\phi$ is not recurrent.

	\begin{corollary}
		\label{recurrentediscreto}
		Let $G$ be  discrete and let $\phi\in  M(A(G)).$ $M_\phi$ is recurrent if and only if $\phi(G)\subseteq \T$.
	\end{corollary}
	
 \begin{remark}
			If $G$ is discrete and countable, since $\phi(G)\subseteq \T$ implies that $\phi(G)$ is an Arbault set, we can find a sequence such that $\lim _kM_\phi^{n_k}(f)=f$ for each $f\in A(G)$ with finite support. This yields that recurrence is equivalent to quasirigidity in the sense of \cite{GLP1,GLP2}.
	\end{remark}
	
	The next result also follows from of \cite[Proposition 2.3]{bulancea25}.
	\begin{proposition}\label{compact-recurrence-ur}
		Let $G$ be a  compact group and let $\phi\in  M(A(G)).$ Then $M_\phi$ is recurrent if and only if $M_\phi$ is uniformly rigid.
	\end{proposition}

\subsubsection{Abelian case}\label{subs:abelian}

We present here a construction based on sums of alternating blocks of Haar measures of a chain of compact subgroups. In every case, the symbols we construct have the same range, an infinite Dirichlet set, and the corresponding multiplication operator on $C_0(G)$ is uniformly rigid. On $A(G)$, however, depending on the parameters, the operator may be uniformly rigid, rigid but not uniformly rigid, or recurrent but not rigid. None of these operators is power bounded.

Throughout this subsection $G$ is a locally compact Abelian group and $\Gamma=\widehat G$. For $\mu\in M(G),$ the convolution operator $C_\mu$ on $L^1(\Gamma)$  and the multiplication operator $M_{\widehat\mu}$ on $A(G)$ are isometrically conjugate. In particular,
\[
 \|M_{\widehat\mu}^{N}-I\|=\|\mu^{*N}-\delta_e\|_{M(\Gamma)}.
\]
Here and below, powers of measures mean convolution powers.


\begin{lemma}\label{haar-block-estimate}
Let $H_1\subsetneq\cdots\subsetneq H_{2r}$ be compact subgroups of an Abelian group, with finite successive indices, and let $\nu_j$ be their normalized Haar measures. Then
\[
 \frac r4\leq\left\|\sum_{j=1}^r(\nu_{2j-1}-\nu_{2j})\right\|\leq 2r.
\]
\end{lemma}
\begin{proof}
  Let $\mu = \sum_{j=1}^r (\nu_{2j-1} - \nu_{2j})$. The upper bound is immediate.  To establish the lower bound, consider the pairwise disjoint sets $E_1 = H_1$ and $E_j = H_j \setminus H_{j-1}$ for $j \ge 2$. Note that $\nu_j(E_k)=0$ if $j<k$.

  We focus on odd $k$'s. That will suffice. To that end, fix $k$, $1\leq k\leq r$.
For each $j$, $k< j \leq r$,  $E_{2k-1}\subseteq H_{2j-1}$. Now, as the restriction of $\nu_{2j}$ to $H_{2j-1}$ is invariant under translations in $H_{2j-1}$, it has to be a multiple of $\nu_{2j-1}$. Checking against $H_{2j-1}$, we see that $\restr{\nu_{2j}}{H_{2j-1}}=\frac{1}{d_{2j-1}}\nu_{2j-1}$.   It follows that  $\bigl(\nu_{2j-1} - \nu_{2j}\bigr)(E_{2k-1})\geq  0$ for every $k<j\leq r$.

If $k=j>1$,
  \begin{align*}
    \bigl(\nu_{2k-1} - \nu_{2k}\bigr)(E_{2k-1})&=  \nu_{2k-1} \left(H_{2k-1}\setminus  H_{2k-2}\right)-\nu_{2k} \left(H_{2k-1}\setminus H_{2k-2}\right)\notag\\
    &=\left(1-\frac{1}{d_{2k-2}}\right)-\left(\dfrac{1}{d_{2k-1}}-\dfrac{1}{d_{2k-1}d_{2k-2}}\right)\\
&=\left(1-\frac{1}{d_{2k-2}}\right)\left(1-\frac{1}{d_{2k-1}}\right)\geq \frac14. \end{align*}
When $k=j=1$, we have that $(\nu_1 - \nu_2)(E_1) = 1 - d_1^{-1} \ge 1/2$.

Combining these cases,  one sees that for each $1\leq k \leq r$, $\mu(E_{2k-1})=\sum_{j=k}^r  (\nu_{2j-1}-\nu_{2j})(E_{2k-1})\geq \dfrac14$. Therefore,
\[
\left\|\mu \right\| \ge \sum_{k=1}^r \mu(E_{2k-1}) \ge \frac{r}{4}.\]
\end{proof}

The construction below follows the same pattern as
Proposition~\ref{leinert}: the symbol equals $1$ outside a sequence of pairwise disjoint sets $(A_n)_n$ and takes the value
$e^{2\pi i/n!}$ on $A_n$. In the Leinert-set construction,
bounded functions supported on such sets have uniformly
controlled multiplier norms. Here, varying the lengths of
the Haar-measure blocks constructed changes the norm control and
produces different dynamical behaviours with the same range.

In what follows, we will use the shortcut $\zeta_n=e^{\frac{2\pi i}{n!}}$.
\begin{proposition}\label{abelian-block-construction}
Assume $G$ is a locally compact Abelian group such  that $\Gamma=\widehat{G}$ contains a strictly increasing chain $(\Gamma_j)_{j\geq1}$ of compact subgroups with finite successive indices. Let $\nu_j$ denote normalized Haar measure on $\Gamma_j$. Choose positive integers $\ell_n$ such that $ \sum_{n=1}^{\infty}\frac{\ell_n}{n!}<\infty$
and set
\begin{equation*}\label{block-measures}
 \sigma_n=\sum_{r=0}^{\ell_n-1}(\nu_{b_n+2r}-\nu_{b_n+2r+1}),
 \qquad \mu=\delta_e+\sum_{n\geq1}(\zeta_n-1)\sigma_n,
\end{equation*}
where  $b_1=1, \ b_{n+1}=b_n+2\ell_n, \ n\in \N.$ Then $\mu\in M(\Gamma)$ and $\phi=\widehat\mu\in B(G)$ satisfies
\[
 \phi(G)=\{e^{2\pi i/n!}:n\geq1\}.
\]
The operator $M_\phi^0$ is uniformly rigid, while $M_\phi$ is not power bounded. Moreover:
\begin{enumerate}
\item[(i)] If $\ell_n=1$ for every $n$, then $M_\phi$ is uniformly rigid.
\item[(ii)] If $\ell_n=n$, then $M_\phi$ is not uniformly rigid. If the connected components of $G$ are open, it is rigid along $(N!)_{N\geq1}$.
\item[(iii)] If $\ell_n=n^2$, then $M_\phi$ is not rigid. If the connected components of $G$ are open, it is recurrent.
\end{enumerate}
\end{proposition}
\begin{proof}
The series defining $\mu$ converges in the  variation norm, since $\|\sigma_n\|\leq2\ell_n$ and $|\zeta_n-1|\leq2\pi/n!$. Put $G_j=\Gamma_j^\perp$. These are strictly decreasing open subgroups of $G$, so $\widehat\nu_j=\Cf_{G_j}\in B(G)$. Hence, for every $n\in \N,$
\begin{equation*}\label{mutrans}
 \widehat\sigma_n=\Cf_{A_n},\ \
 A_n=\bigcup_{r=0}^{\ell_n-1}(G_{b_n+2r}\setminus G_{b_n+2r+1}), \ \ \text{and }\ \phi=1+\sum_{n\geq1}(\zeta_n-1)\Cf_{A_n}.
\end{equation*}
The sets $A_n$ are nonempty and pairwise disjoint. As $e\notin\bigcup_n A_n$,  $1$ belongs to the range. This proves the assertion about $\phi(G),$ which yields uniform rigidity on $C_0(G)$ by Proposition \ref{nonunitalcarac}(iii). The infinite range and Corollary \ref{paffdense} show that $M_\phi$ is not power bounded.

For every fixed positive integer $q$, injectivity of the Fourier--Stieltjes transform gives
\begin{equation}\label{block-powers}
 \mu^{*q}-\delta_e=\sum_{n\geq1}(\zeta_n^q-1)\sigma_n.
\end{equation}
This series also converges in norm, as $|\zeta_n^q-1|\leq2\pi q/n!$.
For a positive integer $q$, put
\[m(q)=\min\{n\geq1:n!\text{ does not divide }q\}.\]
Then $m(q) \ge 2$ and $(m(q)-1)!$ divides $q$, so we can write $q = k (m(q)-1)!$ for some positive integer $k$ with $k \not\equiv 0 \pmod{m(q)}$. This implies that
\[
\zeta_{m(q)}^q = \exp\left( \frac{2\pi i k (m(q)-1)!}{m(q)!} \right) = \exp\left( \frac{2\pi i k}{m(q)} \right).
\]
 Thus,
\[
|\zeta_{m(q)}^q - 1| = \left| \exp\left(\frac{2\pi i k}{m(q)}\right) - 1 \right| \ge \left| \exp\left(\frac{2\pi i}{m(q)}\right) - 1 \right| = 2\sin\left(\frac{\pi}{m(q)}\right).
\]
Since $\sin(x) \ge \frac{2}{\pi}x$ for all $x \in [0, \pi/2]$ and $\frac{\pi}{m(q)} \le \frac{\pi}{2}$, we obtain the   lower bound
\[
|\zeta_{m(q)}^q - 1| \ge  \frac{4}{m(q)}.
\]

To isolate the $m(q)$-th component in \eqref{block-powers}, we take into account that $\nu_i * \nu_j = \nu_{\max(i,j)}$ for all $i, j \ge 1$. Consider  $\eta_q = \delta_e - \nu_{b_{m(q)+1}} \in M(\Gamma)$. As $\sigma_n = \sum_{r=0}^{\ell_n-1} (\nu_{b_n+2r} - \nu_{b_n+2r+1})$, for $n > m(q)$ the measures $\nu_j$ defining $\sigma_n$ satisfy $j \ge b_{n} \ge b_{m(q)+1}$. Consequently, for all $n > m(q)$, the idempotency relation yields
\[
\sigma_n * \nu_{b_{m(q)+1}} = \sigma_n,
\]
which implies that $\sigma_n * \eta_q = \sigma_n * (\delta_e - \nu_{b_{m(q)+1}}) = 0$. For $n = m(q)$, the indices appearing in $\sigma_{m(q)}$ satisfy $b_{m(q)}+2r < b_{m(q)}+2r+1 \le b_{m(q)}+2\ell_{m(q)} = b_{m(q)+1}$. Hence,
\[
\sigma_{m(q)} * \nu_{b_{m(q)+1}} = \sum_{r=0}^{\ell_{m(q)}-1} (\nu_{b_{m(q)+1}} - \nu_{b_{m(q)+1}}) = 0.
\]
Consequently, $\sigma_{m(q)} * \eta_q = \sigma_{m(q)} * (\delta_e - \nu_{b_{m(q)+1}}) = \sigma_{m(q)}$. Since $\zeta_n^q = 1$ for all $n < m(q)$, all earlier terms in \eqref{block-powers} vanish identically. Therefore, convolving \eqref{block-powers} with $\eta_q$ only preserves the $m(q)$-th block term, yielding
\[
(\mu^{*q} - \delta_e) * \eta_q = (\zeta_{m(q)}^q - 1) \sigma_{m(q)}.
\]
 As $\|\eta_q\| \le   2$,  Lemma~\ref{haar-block-estimate} yields
\[
2 \|\mu^{*q} - \delta_e\| \ge \|\mu^{*q} - \delta_e\| \|\eta_q\| \ge \|(\mu^{*q} - \delta_e) * \eta_q\| \ge \frac{4}{m(q)} \|\sigma_{m(q)}\| \ge \frac{\ell_{m(q)}}{m(q)},
\]
which  yields the  lower bound
\begin{equation}\label{block-lower-bound}
	\|M_\phi^q - I\|= \|\mu^{*q} - \delta_e\| \ge \frac{\ell_{m(q)}}{2m(q)}.
\end{equation}
In order to prove   (i), we have
\[
 \|\mu^{*N!}-\delta_e\|\leq4\pi\sum_{n>N}\frac{N!}{n!}
 \leq\frac{4\pi}{N+1}\sum_{j\geq0}\frac1{(N+2)^j}\longrightarrow0.
\]
For  (ii), \eqref{block-lower-bound} gives $\|M_\phi^q-I\|\geq1/2$ for every $q$. On the other hand,
\[
 \|\mu^{*N!}-\delta_e\|\leq4\pi\sum_{n>N}\frac{nN!}{n!}
 \leq4\pi\sum_{j\geq0}\frac1{(N+1)^j}\leq8\pi.
\]
Thus $(M_\phi^{N!})$ is bounded. If the connected components of $G$ are open, as $\phi(G)$ is Dirichlet, Remark~\ref{tauberianquasirigid} gives $M_\phi^{N!}u\to u$ for every compactly supported $u\in A(G)$, which form a dense subset of $A(G),$  and boundedness extends convergence to all $u\in A(G)$.

For (iii), assume that $M_\phi$ is rigid along a sequence $(q_k)_k$.
Fix $n \ge 1$ and choose $x \in A_n$. There exists  $u \in A(G)$ such that $u(x) \ne 0$. Thus, we obtain
\[
 \left| (\phi(x)^{q_k} - 1) u(x) \right| \le \|\phi^{q_k} u - u\|_\infty \le \|\phi^{q_k} u - u\|_{A(G)} \xrightarrow[k \to \infty]{} 0.
\]
As $\phi(x) = \zeta_n = e^{2\pi i / n!}$, we get $\lim_{k \to \infty} \zeta_n^{q_k} = 1.$
Because $\zeta_n$ is an $n!$-th root of unity, there exists $k_0\in \N$ such that $\zeta_n^{q_k} = 1$,  i.e., $n!$ divides $q_k$, for all   $k\geq k_0$.  This implies $m(q_k) > n $ for all $k \ge k_0$. Since $n$ was arbitrary, it follows that $\lim_{k \to \infty} m(q_k) = \infty$. Substituting $\ell_n = n^2$ into the lower bound \eqref{block-lower-bound} yields
\[
\|M_\phi^{q_k} - I\| \ge \frac{\ell_{m(q_k)}}{2 m(q_k)} = \frac{m(q_k)}{2} \xrightarrow[k \to \infty]{} \infty,
\]
contradicting the Uniform Boundedness Principle, as $M_\phi^{q_k}\rightarrow I$ in the strong operator topology.

Finally, as $\phi$ is unimodular and  constant on each connected component, if the connected components of $G$ are open, Theorem~\ref{rec_equiv_locA} implies that $M_\phi$ is recurrent.
\end{proof}


We now explain under which circumstances  the chain required in Proposition~\ref{abelian-block-construction} is available. First, we show it always exists when the connected component   of the identity  is not open.

\begin{lemma}\label{lemma:chain_construction}
If the connected component $C_0$ of the identity of $G$ is not open, then $\Gamma=\widehat{G}$ contains a strictly increasing sequence $(\Gamma_j)_{j\geq1}$ of compact subgroups such that $[\Gamma_{j+1}:\Gamma_j]<\infty$ for every $j\in \N$.
\end{lemma}
\begin{proof}
By the Principal  Structure Theorem for locally compact Abelian groups, $G\cong\R^d\times H$, where $H$ has a compact open subgroup $K$ (see    \cite[Theorem 23.25]{HewittRoss} or \cite[Corollary 3]{Morris}). By  hypothesis,  the connected component of the identity of $H$, $C_0^H,$ is not open. Since $K$ is a closed an open subgroup, we get $C_0^H\subseteq K,$ thus  the identity component of $K$, $C_0^K,$  coincides with that of $H$. As $C_0^K$ is not open,  $K/C_0^K$ is an infinite compact totally disconnected Abelian group, and its dual is an infinite discrete torsion group. Choose a strictly increasing chain $(F_j)$ of finite subgroups of this discrete group $\widehat{K/C_0^K}$. We regard these subgroups $F_j$ as subgroups of  $\widehat K$ that are in $ (C_0^K)^{\perp}$, that is they take the constant value 1 on $C_0^K$.  Since $K$ is an open subgroup, the restriction map $R:\widehat H\to\widehat K$ has compact open kernel $K^\perp$. Consequently the subgroups $R^{-1}(F_j)$, being  a finite union of cosets of $K^\perp$,
 are compact  subgroups of $\widehat H$ with
\[
 [R^{-1}(F_{j+1}):R^{-1}(F_j)]=[F_{j+1}:F_j].
\]
Embedding these subgroups into $\widehat{G} \cong \mathbb{R}^d \times \widehat{H}$ via $\Gamma_j = \{0\} \times R^{-1}(F_j)$ yields the desired chain of compact subgroups in $\Gamma$.
\end{proof}

\begin{example}\label{dyadic-chain-examples}
	{\em Let
\[
 K=\prod_{j\geq1}\Z/2\Z
\]
with the product topology, be the compact dyadic Cantor group. Its dual is the discrete group $D=\bigoplus_{j\geq1}\Z/2\Z$. In $D$, take $\Gamma_j$  the subgroup generated by the first $j$ coordinate vectors. The successive indices are two. The construction in Proposition \ref{abelian-block-construction} applies to $G=K$, whose identity component $\{0\}$ is not open. For $\ell_n=1$ it gives a uniformly rigid multiplier with infinite Dirichlet range. For $\ell_n=n$ and $\ell_n=n^2$, the resulting multipliers $M_\phi$  are not recurrent, since they are not uniformly rigid, and on the unital algebra $A(K)$ recurrence and uniform rigidity are equivalent by Proposition~\ref{compact-recurrence-ur}. In every case $M_\phi^0$ is uniformly rigid.}
\end{example}

The existence of a chain of compact subgroups as required by Proposition~\ref{abelian-block-construction} does not actually need the hypothesis on the connected components used in Lemma~\ref{lemma:chain_construction}. We now show that such chains may also exist when the identity component is open. The resulting examples separate the three dynamical properties under study and show that the equivalence between recurrence and rigidity for power bounded multipliers on algebras with $\sigma$-compact spectrum (Corollary~\ref{sigmacompact}) fails without the power-boundedness assumption, even in the case of the countable discrete spectrum $\mathbb{Z}$.

\begin{theorem}
	 Let $G$ be a locally compact Abelian group with open connected components and that $\widehat G$ contains a chain as in Proposition~\ref{abelian-block-construction}. Then there are non-power-bounded multipliers $M_\phi$, $\phi\in B(G)$  with $M_\phi^0$ uniformly rigid such that:
\begin{itemize}
  \item[(i)] $M_\phi$ is rigid but not uniformly rigid.
  \item[(ii)] $M_\phi$ is recurrent but not rigid
\end{itemize}   This applies in particular to $G=\Z$, $G=\R^d\times\Z$, and $G=\bigoplus_{j\geq1}\Z/2\Z$ with the discrete topology.
\end{theorem}

\begin{proof}
The assertion follows by Proposition~\ref{abelian-block-construction}, taking $\ell_n=n$ and $\ell_n=n^2$, respectively. Let us verify that $G=\Z$, $G=\R^d\times\Z$, and $G=\bigoplus_{j\geq1}\Z/2\Z$ admit the required chain in their dual.\\
For $G=\Z$, the dual group is $\T$. Take
\[
\Gamma_j=\{z\in\T:z^{2^j}=1\}.
\]
These are compact finite subgroups with $[\Gamma_{j+1}:\Gamma_j]=2$. All three conclusions of Proposition~\ref{abelian-block-construction} apply, since $G$ is discrete. More generally, they apply to $G=\R^d\times\Z$, whose connected components are open, by taking $\{0\}\times\Gamma_j$ in the dual.

Consider now $G=\bigoplus_{j\geq1}\Z/2\Z$ with the discrete topology, so that  $\widehat{G}=K$ in Example \ref{dyadic-chain-examples}. Define the finite (compact) subgroups
\[
\Gamma_j=\{x\in K:x_k=0\text{ for every }k>j\}.
\]
Their successive indices are two. Although these subgroups are not open in $K$, each is open in the next one, which is precisely the finite-index condition used in the construction. Since $G$ is discrete, the cases $\ell_n=1,n,n^2$ yield, respectively, uniform rigidity, rigidity without uniform rigidity, and recurrence without rigidity.
\end{proof}
Theorem \ref{rec_equiv_locA} shows that, when the connected components are open, $M_\phi$ is recurrent if and only if $M_\phi^0$ is. This equivalence fails in every group whose connected components are not open.

\begin{theorem}\label{C0_open_equiv_rec}
Let $G$ be a locally compact Abelian group. The following are equivalent:
\begin{enumerate}
\item[(i)] There exists $\phi\in B(G)$ such that $M_\phi^0$ is uniformly rigid but $M_\phi$ is not recurrent.
\item[(ii)] The connected component of the identity of $G$ is not open.
\end{enumerate}
\end{theorem}
\begin{proof}
If the identity component is open, Theorem~\ref{rec_equiv_locA} excludes (i). Conversely, write $G\cong\R^d\times H$, where $H$ has a compact open subgroup $K$. Then $C_0^K$ is not open in $K$. Apply Lemma~\ref{lemma:chain_construction} and Proposition~\ref{abelian-block-construction}, with $\ell_n=n$, to $K$. This gives $\widetilde\phi\in B(K)=A(K)$ for which $M_{\widetilde\phi}^0$ is uniformly rigid but $M_{\widetilde\phi}$ is not uniformly rigid. Since $A(K)$ is unital, Theorem~\ref{compact-recurrence-ur} shows that $M_{\widetilde\phi}$ is not recurrent.
The extension $f$ of $\widetilde\phi$ by zero from the open subgroup $K$ to $H$ belongs to $B(H)$ \cite[Proposition 7.1.2]{kanilaubook}. Set $\phi_H=f+1-\Cf_K$ and $\phi(x,h)=\phi_H(h)$. Then $\phi\in B(G)$ and its range is the infinite Dirichlet set $\{\zeta_n:n\geq1\}$, so $M_\phi^0$ is uniformly rigid. The restriction map $A(G)\to A(\{0\}\times K)\cong A(K)$ is continuous and surjective \cite[Corollary 2.6.5]{kanilaubook} and intertwines $M_\phi$ with $M_{\widetilde\phi}$. Recurrence of $M_\phi$ would therefore imply recurrence of $M_{\widetilde\phi}$, a contradiction.
\end{proof}

Sometimes every uniformly rigid multiplier of  $A(G)$  is power bounded.
\begin{proposition}
	Let $G$ be a compact  Abelian group. The following are equivalent:
	\begin{itemize}
		\item[(i)] The connected component of the identity of $G$ is not open.
		\item[(ii)] There exists $\phi \in B(G)$ not power-bounded such that $M_\phi : A(G) \to A(G)$ is uniformly rigid.
	\end{itemize}
\end{proposition}

\begin{proof}
	(i) $\Rightarrow$ (ii). If the connected component of the identity  $C_0$ of $G$ is not open, Lemma~\ref{lemma:chain_construction} provides a sequence $(\Gamma_n)_{n \ge 1} \subseteq \widehat{G}$ of compact subgroups with $[\Gamma_{n+1} : \Gamma_n] < \infty$. Applying Proposition~\ref{abelian-block-construction} with $\ell_n = 1$ yields a non-power-bounded multiplier $\phi \in B(G)$ such that $M_\phi$ is uniformly rigid.\\
		(ii) $\Rightarrow$  (i). If $C_0$ is open, then $G = \bigcup_{g \in F} g C_0$ for some finite set $F \subset G$. If $M_\phi$ is uniformly rigid, then $\phi(G)$ is a Dirichlet subset of $\T.$ Since $C_0$ is connected, each $\phi(g C_0)$ is constant. Thus, $\phi(G)$ is finite, so $\phi$ is power-bounded by Proposition \ref{rangofinito}.
\end{proof}

 These results have a concrete interpretation for convolution operators $C_\mu f=\mu*f$ on $L^1(\T)$, identified by the Fourier transform with multipliers of $A(\Z)$.

 \begin{corollary}\label{circle-convolution}
 	Let $\mu\in M(\mathbb{T})$ and consider    $C_\mu\colon L^1(\mathbb{T})\to L^1(\mathbb{T}),$  $C_\mu f=\mu*f$.
 	\begin{enumerate}
 		\item[(i)] $C_\mu$ is recurrent if and only if $|\widehat\mu(k)|=1$ for every $k\in\Z$.
 		\item[(ii)] If $C_\mu$ is power bounded, it is recurrent if and only if it is rigid, and it is uniformly rigid if and only if $\widehat\mu(\Z)$ is a finite subset of $\T$.
 		\item[(iii)] There exist non-power-bounded convolution operators on $L^1(\T)$ which are uniformly rigid, rigid but not uniformly rigid, and recurrent but not rigid. All three can be chosen with the same Fourier--Stieltjes range $\widehat\mu(\Z)=\{e^{2\pi i/n!}:n\geq1\}$; their associated operators on $c_0(\Z)$ are uniformly rigid.
 	\end{enumerate}
 \end{corollary}

\begin{proof}
The Fourier transform identifies $L^1(\T)$ isometrically with $A(\Z)$. Apply Corollary~\ref{recurrentediscreto}, Corollary~\ref{sigmacompact}, Theorem~\ref{amenableur}, and Proposition~\ref{abelian-block-construction}, using the chain of subgroups of $\T$ of orders $2^j$, $j\in \N.$
\end{proof}
In the literature, a measure on a locally compact Abelian group whose Fourier–Stieltjes transform takes values in $\T$ has been  called a \emph{unitary measure}. For a discussion of their properties, see \cite{corw74} and the references therein.

\section*{Use of artificial intelligence}\vspace{-4mm}
AI-assisted tools (ChatGPT, OpenAI, GPT-6) suggested the initial approach to Proposition ~\ref{leinert}. Beyond this result, AI tools were used only to check computations and to assist with editing. The authors have    independently verified all content and takes full responsibility for its accuracy.
\vspace{3mm}\par

\noindent{\bf Acknowledgements}. The research of the first and the third  authors was partially supported by the projectPID2024-162128NB-I00 funded by MCIN/AEI/10.13039/501100011033 and by “ERFD Away of making Europe” and GVA AICO/2021/17.	Research of the second author  was partially supported by UJI, GACUJIMB/2025/04.

\end{document}